\documentclass{article}
\usepackage{amsmath,amsthm,amssymb}
\usepackage{times}
\usepackage{enumerate}
\usepackage{amssymb}
\usepackage{graphicx}
\usepackage{color}
\usepackage{cite}

\makeatletter
\def\author#1{\gdef\autrun{\def\and{\unskip, }#1}\gdef\@author{#1}}

\makeatother

\def\keywords#1{\par\medskip
\noindent\textbf{Keywords.} #1}
\def\subjclass#1{\par\smallskip
\noindent\textbf{MSC (2010):} #1}

\newtheorem{thm}{Theorem}[section]

\newtheorem{lem}[thm]{Lemma}

\newtheorem{prop}[thm]{Proposition}

\theoremstyle{definition}

\newtheorem{rem}[thm]{Remark}
\newtheorem{exa}[thm]{Example}

\numberwithin{equation}{section}

\makeatletter
\def\author#1{\gdef\autrun{\def\and{\unskip, }#1}\gdef\@author{#1}}

\makeatother

\makeatletter
\let\@fnsymbol\@alph
\makeatother

\newtheorem*{notations}{Notations}

\DeclareMathOperator*{\esssup}{ess\,sup}
\DeclareMathOperator*{\essinf}{ess\,inf}

\begin{document}

\title{Radial solutions for the bilaplacian equation with vanishing or singular radial potentials}

\author{Marino Badiale\thanks{Dipartimento di Matematica ``Giuseppe Peano'', Universit\`{a} degli Studi di
Torino, Via Carlo Alberto 10, 10123 Torino, Italy. 
e-mail: \texttt{marino.badiale@unito.it}}
\textsuperscript{,}\thanks{Partially supported by the MIUR 2015 PRIN project ``Variational methods, with applications to problems in mathematical physics and geometry''.}
\ -\ Stefano Greco
\ -\ Sergio Rolando\thanks{Dipartimento di Matematica e Applicazioni, Universit\`{a} di Milano-Bicocca,
Via Roberto Cozzi 53, 20125 Milano, Italy. e-mail: \texttt{sergio.rolando@unito.it}. Member of the Gruppo Nazionale di Alta Matematica (INdAM).}%
\ \textsuperscript{,\,b}
}

\date{
}
\maketitle

\begin{abstract}
\noindent Given three measurable functions $V\left(r \right)\geq 0$, $K\left(r\right)> 0$ and $Q\left(r \right)\geq 0$, $r>0$,
we consider the bilaplacian equation
\[
\Delta^2 u+V(|x|)u=K(|x|)f(u)+Q(|x|) \quad \text{in }\mathbb{R}^N
\]
and we find radial solutions thanks to compact embeddings of radial spaces of Sobolev functions into sum of weighted Lebesgue spaces.

\end{abstract}

\section{Introduction}

This paper is concerned with the following bilaplacian equation
\begin{equation} \label{prob1}
\Delta^2 u+V(|x|)u=K(|x|)f(u)+Q(|x|) \quad \text{in } \mathbb{R}^N
\end{equation}
where $N\geq5$, $\Delta^2 u = \Delta(\Delta u)$ is the $bilaplacian$ operator, the forcing term $Q\geq 0$ and the potentials $V\geq 0$ and $K>0$ satisfy suitable hypotheses, and $f:\mathbb{R}\rightarrow\mathbb{R}$ is a continuous function such that $f(0)=0$. 
In particular we are interested in allowing the potential $V$ to be singular at the origin and/or vanishing at infinity.

Bilaplacian equations arise in describing different physical phenomena, such as the propagation of laser beams in Kerr media or nonlinear oscillations in suspension bridges (see some references in \cite{sun-ecc, bastos-ecc}), and have been extensively studied in the last decades (see e.g. \cite{sun-ecc, sun-ecc-2,bastos-ecc, liang-ecc, carr, demarque-miya} and the references therein). 
In spite of that, equations of type (\ref{prob1}), namely with radial potentials possibly singular at the origin and vanishing at infinity, has been treated only in \cite{carr, demarque-miya} (at least to our knowledge), where the authors essentially consider power type potentials.

For problem (\ref{prob1}) we will obtain several kinds of existence results of radial solutions.
The main technical device for our results is given by some new theorems about compact embeddings of suitable Sobolev spaces into sum of weighted Lebesgue spaces. 
The natural approach in studying Eq. \eqref{prob1} is variational, since
its weak solutions are (at least formally) critical points of a suitable Euler functional, as we will see. 
Then the problem of existence is easily solved if $V$ does not vanish at infinity and $K$ is
bounded, because standard embeddings theorems are available. As we will let $V$ and $K$ to vanish, 
or to go to infinity, as $|x| \rightarrow 0$ and $|x|\rightarrow +\infty$, the usual embeddings theorems for Sobolev spaces are not 
available anymore, and new embedding theorems must be proved. This kind of work has been started in \cite{carr, demarque-miya} (for the bilaplacian equation) and we continue it here, using some new ideas that has been introduced in \cite{bad1,bad2,bad3}.

The main novelty of our approach is two-folded. Firstly, we look for embeddings not into a Lebesgue space 
but into a sum of Lebesgue spaces. 
This allows us to study separately the behavior of $V$ and $K$ at $0$ and $\infty$, and to assume different set of hypotheses about these behaviors. Secondly, we assume hypotheses not on $V$ and $K$ separately but on their ratio, so allowing general kind of asymptotic behavior for the two potentials.

Thanks to this second novelty we obtain embedding results, and thus existence results for Eq. \eqref{prob1}, which extend the ones of \cite{carr, demarque-miya} to more general kinds of potentials. Moreover, thanks to the first novelty, we get new results also for power type potentials (cf. Example \ref{EX} below).

The paper is organized as follows. In Sections \ref{SEC:MAIN} and \ref{SEC: ex} we give our results on compact embeddings and existence of solutions to Eq. \eqref{prob1} respectively. The former will be proved in the Sections \ref{thm1}-\ref{thms45}, the latter in Section \ref{SEC: app pf}.

\begin{notations}
We end this introductory section by collecting
some notations used in the paper.

\noindent $\bullet $ We denote $\mathbb{R}_+ :=(0,+\infty)$ and $\mathbb{R}_- :=(-\infty,0)$.

\noindent $\bullet $ For every $R>0$, we set $B_{R}:=\left\{ x\in \mathbb{R}%
^{N}:\left| x\right| <r\right\} $.

\noindent $\bullet $ For any subset $A\subseteq \mathbb{R}^{N}$, we denote $%
A^{c}:=\mathbb{R}^{N}\setminus A$. 

\noindent $\bullet $ By $\rightarrow $ and $\rightharpoonup $ we
respectively mean \emph{strong} and \emph{weak }convergence.

\noindent $\bullet $ $C_c^{\infty }(\Omega )$ is the space of the
infinitely differentiable real functions with compact support in the open
set $\Omega \subseteq \mathbb{R}^{N}$.

\noindent $\bullet $ If $1\leq p\leq \infty $ then $L^{p}(A)$ and $L_{%
\mathrm{loc}}^{p}(A)$ are the usual real Lebesgue spaces (for any measurable
set $A\subseteq \mathbb{R}^{N}$). If $\rho :A\rightarrow \mathbb{R}_+$ 
is a measurable function, then $L^{p}(A,\rho \left( z\right) dz)$
is the real Lebesgue space with respect to the measure $\rho \left( z\right)
dz$ ($dz$ stands for the Lebesgue measure on $\mathbb{R}^{N}$).

\noindent $\bullet $ $p^{\prime }:=p/(p-1)$ is the H\"{o}lder-conjugate exponent of $p.$

\end{notations}

\section{Main results} \label{SEC:MAIN}

In this section we state our main results on compact embeddings, that we will prove in the following Sections \ref{thm1}-\ref{thms45}. Firstly, we introduce some basic concepts and results. Assume $N\geq 5$ and define $2^{**}:=\frac{2N}{N-4}$. 

By usual Sobolev embeddings, there exists a suitable constant $C>0$ such that 
for all $u \in C_c^\infty(\mathbb{R}^N)$ one has
$$\left\|u\right\|_{L^{2^{**}}} \leq C \left\| D^2 u \right\|_{L^2}
$$
where
\begin{equation} \label{normasolita}
\Vert D^2 u \Vert_{L^2}:=\left( \sum_{|\alpha|=2} \left \| D^\alpha u \right \|^2_{L^2} \right)^{1/2}.
\end{equation}

A basic space to work with is 
\[
D^{2,2}(\mathbb{R}^N):=\left\{u \in L^{2^{**}}(\mathbb{R}^N) : \left\| D^2 u \right\|_{L^2}< +\infty\right\}. 
\]
It is the closure of $C_c^\infty(\mathbb{R}^N)$ in $L^{2^{**}}(\mathbb{R}^N) $ with respect to the norm $\left\| D^2 u \right\|_{L^2}$
and, endowed with such a norm, it is an Hilbert space. 
The bilinear form
$$(u,v)\mapsto \int_{\mathbb{R}^N} \Delta u \Delta v \, dx$$
defines a scalar product on $D^{2,2}(\mathbb{R}^N)$
and the associated norm, that is $\Vert u \Vert_{D^{2,2}}:=\Vert \Delta u \Vert_{L^2}$, is equivalent to \eqref{normasolita}
(see for example \cite{gazz}).
Hence, one can also define $D^{2,2}(\mathbb{R}^N)$ as the closure of $C_c^\infty(\mathbb{R}^N)$ in $L^{2^{**}}(\mathbb{R}^N)$ with respect to the norm $\Vert \Delta u \Vert_{L^2}$. We will be particularly interested in the subspace of radial functions, i.e., 
$$
D^{2,2}_r := D^{2,2}_r(\mathbb{R}^N):=\left\{ u\in D^{2,2}(\mathbb{R}^N): u(x)=u(|x|) \right\},
$$
for which the pointwise estimates given in the following lemma hold (see \cite{nouss} for a proof).

\begin{lem} \label{stimepunt}
For every $u\in D^{2,2}_r(\mathbb{R}^N)$ we have
\begin{equation} \label{stima1}
|u(x)| \leq \frac{2}{N-4} \frac{1}{\sqrt{N\sigma_N}}
\dfrac{\left \| \Delta u \right \|_{L^2}}{|x|^\frac{N-4}{2}} \quad \text{almost everywhere in }\mathbb{R}^N
\end{equation}
and
\begin{equation} \label{stima2}
|\nabla u(x)| \leq \frac{1}{\sqrt{N\sigma_N}}
\dfrac{\left \| \Delta u \right \|_{L^2}}{|x|^\frac{N-2}{2}}  \quad \text{almost everywhere in }\mathbb{R}^N
\end{equation}
where $\sigma_N$ denotes the $(N-1)$-dimensional measure of the unit sphere of $\mathbb{R}^N$.
\end{lem}

For any measurable function $V: \mathbb{R}_+\rightarrow [0,+\infty]$, we define the space
$$H^2_V(\mathbb{R}^N):=\left\lbrace u\in D^{2,2}(\mathbb{R}^N): \int_{\mathbb{R}^N} V(|x|)|u|^2 \, dx <\infty \right\rbrace
$$
This is an Hilbert space with scalar product
\begin{equation}\label{inn.prod}
(u,v):=\int_{\mathbb{R}^N} \Delta u \Delta v \, dx +\int_{\mathbb{R}^N} V(|x|)uv \, dx
\end{equation}
and associated norm
$$\left\|  u \right\|:=\left( \int_{\mathbb{R}^N} |\Delta u|^2 +V(|x|)|u|^2 \, dx \right)^\frac{1}{2}.$$
We are interested in finding solutions of \eqref{prob1} in the radial subspace of $H^2_V(\mathbb{R}^N)$, i.e.,
$$H^2_{V,r}=H^2_{V,r}(\mathbb{R}^N):=\left\{ u\in H^2_V(\mathbb{R}^N): u(x)=u(|x|) \right\}.$$

\begin{rem}
By the Sobolev embedding, there is a constant $S_N>0$ such that
\begin{equation} \label{immsobo}
\forall u \in H^2_{V} (\mathbb{R}^N), \quad \left \| u \right \|_{L^{2^{**}}}\leq S_N \left \| u \right \|.
\end{equation}
\end{rem}

\begin{rem}
From the continuous embedding $H^2_{V,r}\hookrightarrow D^{2,2}_{r}(\mathbb{R}^N)$ and inequality \eqref{stima1}, we deduce that there exists a constant $C_N>0$ such that
\begin{equation} \label{stimauso}
\forall u\in H^2_{V,r}(\mathbb{R}^N), \quad |u(x)| \leq C_N \dfrac{\left \|  u \right \|}{|x|^\frac{N-4}{2}} 
\quad \text{almost everywhere in }\mathbb{R}^N.
\end{equation}
\end{rem}

We now introduce the sum of Lebesgue spaces. 
For any measurable function $K: \mathbb{R}_+\rightarrow \mathbb{R}_+$ and $1< q_1\leq q_2 <\infty$, we define
$$L^{q_1}_K+L^{q_2}_K:=L^{q_1}_K(\mathbb{R}^N)+L^{q_2}_K(\mathbb{R}^N):=\left\lbrace u_1+u_2:u_i \in L^{q_i}_K(\mathbb{R}^N), i=1,2 \right\rbrace.$$
This space coincides with the set of measurable functions  $u:\mathbb{R}^N\rightarrow \mathbb{R}$ for which there exists a measurable set $E\subseteq \mathbb{R}^N$ such that $u\in L^{q_1}_K(E)\cap L^{q_2}_K(E^c)$ (where $L^{q_1}_K(E):=L^{q_1}(E,K(|x|)dx)$ and $L^{q_2}_K(E^c):=L^{q_2}(E^c,K(|x|)dx)$)
and it is a Banach space when endowed with the norm 
$$\left \|  u \right \|_{L^{q_1}_K+L^{q_2}_K}:=\inf_{u_1+u_2=u} \max\left\lbrace \left \|  u_1 \right \|_{L^{q_1}_K},\left \|  u_2 \right \|_{L^{q_2}_K} \right\rbrace$$
(see \cite{bad4}). Note that $L^{q}_K$ is continuously embedded into $L^{q_1}_K+L^{q_2}_K$ for all $q\in \left[q_1,q_2 \right]$. 
\smallskip

Our first result is Theorem \ref{teo1} below. It provides sufficient condition to the embeddings we are interested in and uses the following assumptions:

\begin{itemize}
\item[$\left( \mathbf{V}\right) $]  
$V: \mathbb{R}_+\rightarrow[0,+\infty]$ is a measurable function such that $V\in L^1_{loc}(\mathbb{R}_+)$

\item[$\left( \mathbf{K}\right) $]  
$K: \mathbb{R}_+\rightarrow\mathbb{R}_+$ is a measurable function such that $K\in L^s_{loc}(\mathbb{R}_+)$ for some  $s>1$

\item[$\left( \mathcal{S}_{q_{1},q_{2}}^{\prime }\right) $] 
$\exists R_{1},R_{2}>0$ such that $\mathcal{S}_{0}\left( q_{1},R_{1}\right) <\infty $
and $\mathcal{S}_{\infty }\left( q_{2},R_{2}\right) <\infty $

\item[$\left( \mathcal{S}_{q_{1},q_{2}}^{\prime\prime }\right) $] 
$\displaystyle\lim_{R \to 0^+} \mathcal{S}_0(q_1,R)=\lim_{R \to \infty}\mathcal{S}_\infty(q_2,R)=0$
\end{itemize}
where $q_1, q_2$ will be specified each time and $\mathcal{S}_{0},\mathcal{S}_{\infty }$ are the functions of $R>0$ and $q>1$ defined as follows:
\begin{equation} \label{S0}
\mathcal{S}_0(q,R):=\sup_{u\in H^2_{V,r}, \, \left \| u \right \|=1} \int_{B_R} K(|x|)|u|^q \, dx,
\end{equation}
\begin{equation} \label{Sinf}
\mathcal{S}_\infty(q,R):=\sup_{u\in H^2_{V,r}, \, \left \| u \right \|=1} \int_{\mathbb{R}^N\setminus B_R} K(|x|)|u|^q \, dx.
\end{equation}
Notice that $\mathcal{S}_0(q,\cdot)$ is increasing, while $\mathcal{S}_\infty(q,\cdot)$ is decreasing.

\begin{thm} \label{teo1}
Assume \textsc{(\textbf{V})} and  \textsc{(\textbf{K})}, and let $q_1,q_2>1$.
\begin{enumerate}
\item If $\left( \mathcal{S}_{q_{1},q_{2}}^{\prime }\right) $ holds, 
then $H^2_{V,r}(\mathbb{R}^N)$ is continuosly embedded into $L^{q_1}_K(\mathbb{R}^N)+L^{q_2}_K(\mathbb{R}^N)$.
\item If $\left( \mathcal{S}_{q_{1},q_{2}}^{\prime\prime }\right) $ holds, 
then $H^2_{V,r}(\mathbb{R}^N)$ is compactly embedded into $L^{q_1}_K(\mathbb{R}^N)+L^{q_2}_K(\mathbb{R}^N)$.

\end{enumerate}
\end{thm}

We now define two new functions of $R>0$ and $q>1$ as follows:
\begin{equation} \label{R0}
\mathcal{R}_0(q,R):=\sup_{u\in H^2_{V,r}, \, h\in H^2_{V,r}, \, \left \| u \right \|=\left \| h \right \|=1} \int_{B_R} K(|x|)|u|^{q-1}|h| \, dx,
\end{equation}
\begin{equation} \label{Rinf}
\mathcal{R}_\infty(q,R):=\sup_{u\in H^2_{V,r}, \, h\in H^2_{V,r}, \, \left \| u \right \|=\left \| h \right \|=1} \int_{\mathbb{R}^N\setminus B_R} K(|x|)|u|^{q-1}|h| \, dx.
\end{equation}
Note that $\mathcal{R}_0(q,\cdot)$ is increasing, while $\mathcal{R}_\infty(q,\cdot)$ is decreasing.
Furthermore, for any $(q,R)$ we have $\mathcal{S}_0(q,R)\leq \mathcal{R}_0(q,R)$ and $\mathcal{S}_\infty(q,R)\leq \mathcal{R}_\infty(q,R)$, so that $\left( \mathcal{S}_{q_{1},q_{2}}^{\prime\prime }\right) $ is a consequence of the following stronger condition:
\begin{itemize}
\item[$\left( \mathcal{R}_{q_{1},q_{2}}^{\prime\prime }\right) $] 
$\displaystyle\lim_{R \to 0^+} \mathcal{R}_0(q_1,R)=\lim_{R \to \infty}\mathcal{R}_\infty(q_2,R)=0$.
\end{itemize}
In our next results we look for concrete conditions ensuring $\left( \mathcal{R}_{q_{1},q_{2}}^{\prime\prime }\right) $ and thus the compactness of the embedding $H^2_{V,r}(\mathbb{R}^N) \hookrightarrow L^{q_1}_K(\mathbb{R}^N)+L^{q_2}_K(\mathbb{R}^N)$.

Our first results in this direction are Theorems \ref{teo2} and \ref{teo3} below. 
For any $\alpha\in \mathbb{R}$ and $\beta \in [0,1]$, define the functions 
$$\alpha^*(\beta):=\max\left\lbrace 4\beta-2-\dfrac{N}{2},-(1-\beta)N \right\rbrace=\begin{cases}4\beta-2-N/2 & \text{if} \ \ 0 \leq \beta \leq 1/2\\
-(1-\beta)N & \text{if} \ \ 1/2 \leq \beta \leq 1 \end{cases}$$
and
$$q^*(\alpha,\beta):=2\dfrac{\alpha-4\beta+N}{N-4}.$$

\begin{thm} \label{teo2}
Assume \textsc{(\textbf{V})} and \textsc{(\textbf{K})}. 
Assume that $\exists R_1>0$ such that $V(r)<+\infty$ for almost every $r\in (0,R_1)$ and
\begin{equation} \label{ipo1T2}
\esssup_{r\in(0,R_1)} \ \dfrac{K(r)}{r^{\alpha_0} V(r)^{\beta_0}}<+\infty 
\quad \text{for some } 0\leq\beta_0\leq1 \text{ and }\alpha_0>\alpha^*(\beta_0).
\end{equation}
Then $\displaystyle\lim_{R\to0^+} \mathcal{R}_0(q_1,R)=0$ for all $q_1\in \mathbb{R}$ such that
\begin{equation} \label{ipo2T2}
\max\left\lbrace 1,2\beta_0 \right\rbrace < q_1<q^*(\alpha_0,\beta_0).
\end{equation}
\end{thm}

\begin{thm} \label{teo3}
Assume \textsc{(\textbf{V})} and \textsc{(\textbf{K})}. 
Assume that $\exists R_2>0$ such that $V(r)<+\infty$ for almost every $r>R_2$ and
\begin{equation} \label{ipo1T3}
\esssup_{r>R_2} \ \dfrac{K(r)}{r^{\alpha_\infty} V(r)^{\beta_\infty}}<+\infty 
\quad \text{for some } 0\leq\beta_\infty\leq1 \text{ and }\alpha_\infty\in \mathbb{R}.
\end{equation}
Then $\displaystyle \lim_{R\to+\infty} \mathcal{R}_\infty(q_2,R)=0$ for all $q_2\in \mathbb{R}$ such that 
\begin{equation} \label{ipo2T3}
q_2>\max\left\lbrace 1,2\beta_\infty,q^*(\alpha_\infty,\beta_\infty) \right\rbrace.
\end{equation}
\end{thm}

Note that for all $(\alpha,\beta)\in \mathbb{R}\times[0,1]$, we have
$$\max\{1,2\beta,q^*(\alpha,\beta)\}=\begin{cases} q^*(\alpha,\beta) & \text{if} \ \ \alpha \geq \alpha^*(\beta)\\
\max\{1,2\beta\} & \text{if} \ \ \alpha \leq \alpha^*(\beta) \end{cases}.$$

\begin{rem} \label{ossT2T3}

 It is easy to check that the inequalities $\max\left\lbrace 1,2\beta_0 \right\rbrace<q^*(\alpha_0,\beta_0)$ and $\alpha_0>\alpha^*(\beta_0)$ are equivalent. Hence, in \eqref{ipo2T2}, the inequality  $\max\left\lbrace 1,2\beta_0 \right\rbrace<q^*(\alpha_0,\beta_0)$ is automatically satisfied.

\end{rem}

In the next two theorems we assume stronger hypotheses than those of Theorems \ref{teo2} and \ref{teo3}, and we get stronger results. For all $\alpha\in \mathbb{R}$, $\beta\leq1$ and $\gamma\in \mathbb{R}$, define
\begin{equation} \label{le2q*}
q_{*}(\alpha,\beta,\gamma):=2\dfrac{\alpha-\gamma\beta+N}{N-\gamma} \quad\text{and}\quad
q_{**}(\alpha,\beta,\gamma):=2\dfrac{2\alpha+(1-2\beta)\gamma+2(N-2)}{2(N-2)-\gamma}.
\end{equation}
Notice that $q_*$ is defined for $\gamma\neq N$, while $q_{**}$ for $\gamma\neq 2(N-2)$.

\begin{thm} \label{teo4}
Assume \textsc{(\textbf{V})} and \textsc{(\textbf{K})}. 
Assume that $\exists R_2>0$ such that $V(r)<+\infty$ for almost every $r>R_2$ and
\begin{equation} \label{esssup4}
\esssup_{r>R_2} \ \dfrac{K(r)}{r^{\alpha_\infty} V(r)^{\beta_\infty}}<+\infty 
\quad \text{for some } 0\leq\beta_\infty\leq1 \text{ and }\alpha_\infty\in \mathbb{R}
\end{equation}
and
\begin{equation} \label{essinf4}
\essinf_{r>R_2} \ r^{\gamma_\infty}V(r)>0 \quad \text{for some } \gamma_\infty\leq 4.
\end{equation}
Then $\displaystyle \lim_{R\to+\infty} \mathcal{R}_\infty(q_2,R)=0$ for all $q_2\in \mathbb{R}$ such that 
\begin{equation} \label{tesi4}
q_2>\max\left\lbrace 1,2\beta_\infty,q_*,q_{**} \right\rbrace,
\end{equation}

\noindent where $q_*=q_*(\alpha_\infty,\beta_\infty,\gamma_\infty)$ and $q_{**}=q_{**}(\alpha_\infty,\beta_\infty,\gamma_\infty)$.
\end{thm}

To give the statement of our last embedding result, we need to define a subset $\mathcal{A}_{\beta,\gamma}$ of the plane $(\alpha,q)$. Recalling the definitions of $q_*=q_*(\alpha,\beta,\gamma)$ and $q_{**}=q_{**}(\alpha,\beta,\gamma)$ in \eqref{le2q*}, we set
\begin{align} \label{descriA}
\mathcal{A}_{\beta,\gamma}&:=\{(\alpha,q): \max\{1,2\beta \}<q<\min\{q_*,q_{**}\} \} & \text{if}& \ 4\leq\gamma<N, \notag \\
\mathcal{A}_{\beta,\gamma}&:=\{(\alpha,q): \max\{1,2\beta \}<q<q_{**},\alpha>-(1-\beta)N \} & \text{if}& \ \gamma=N, \notag \\
\mathcal{A}_{\beta,\gamma}&:=\{(\alpha,q): \max\{1,2\beta,q_* \}<q<q_{**} \} & \text{if}& \ N<\gamma<2N-4, \notag \\
\mathcal{A}_{\beta,\gamma}&:=\{(\alpha,q): \max\{1,2\beta,q_* \}<q,\alpha>-(1-\beta)\gamma \} & \text{if}& \ \gamma=2N-4, \notag \\
\mathcal{A}_{\beta,\gamma}&:=\{(\alpha,q): \max\{1,2\beta,q_*,q_{**} \}<q \} & \text{if}& \ \gamma>2N-4. 
\end{align}

\begin{thm} \label{teo5}
Assume \textsc{(\textbf{V})} and \textsc{(\textbf{K})}. 
Assume that $\exists R_1>0$ such that $V(r)<+\infty$ for almost every $r\in (0,R_1)$ and
\begin{equation} \label{esssup5}
\esssup_{r\in(0,R_1)} \ \dfrac{K(r)}{r^{\alpha_0} V(r)^{\beta_0}}<+\infty 
\quad \text{for some } 0\leq\beta_0\leq1 \text{ and }\alpha_0\in\ \mathbb{R}
\end{equation}
and
Then $\displaystyle\lim_{R\to0^+} \mathcal{R}_0(q_1,R)=0$ for all $q_1\in \mathbb{R}$ such that
\begin{equation} \label{stareinA}
(\alpha_0,q_1)\in \mathcal{A}_{\beta_0,\gamma_0}.
\end{equation}
\end{thm}

We end this section with an example that might clarify how to use our
results and compare them with the ones of \cite{carr, demarque-miya}. 
Many other examples can be easily obtained by adapting the ones given in \cite[Section 3]{bad1}.

\begin{exa}   \label{EX}
Consider the potentials 
\[
V\left( r\right) =\frac{1}{r^{a}},\quad K\left( r\right) =\frac{1}{r^{a-1}}%
,\quad a\leq 4.
\]
Since $V$ satisfies (\ref{essinf4}) with $\gamma _{\infty }=a$ , we apply Theorem \ref{teo1} together with
Theorems \ref{teo2} and \ref{teo4}. Assumptions (\ref{ipo1T2}) and (\ref{esssup4}) hold if and only if $\alpha _{0}\leq a\beta _{0}-a+1$ and $%
\alpha _{\infty }\geq a\beta _{\infty }-a+1$. According to (\ref{ipo2T2}) and (%
\ref{tesi4}), it is convenient to choose $\alpha _{0}$ as large as possible
and $\alpha _{\infty }$ as small as possible, so we take 
\[
\alpha _{0}=a\beta _{0}-a+1,\quad \alpha _{\infty }=a\beta _{\infty }-a+1.
\]
Then $q^{*}=q^{*}\left( \alpha _{0},\beta _{0}\right) $, $q_{*}=q_{*}\left(
\alpha _{\infty },\beta _{\infty },a\right) $ and $q_{**}=q_{**}\left(
\alpha _{\infty },\beta _{\infty },a\right) $ are given by 
\begin{equation}
q^{*}=2\frac{N-a+1-\left( 4-a\right) \beta _{0}}{N-4},\quad 
q_{*}=2\frac{N-a+1}{N-a}\quad \textrm{and}\quad 
q_{**}=2\frac{2N-a-2}{2N-a-4},  \label{qqq}
\end{equation}
where $a\leq 4$ implies $q_{*}\leq q_{**}$. Note that $\alpha _{0}>\alpha
^{*}\left( \beta _{0}\right) $ for every $\beta _{0}$. Since $q^{*}$ is
decreasing in $\beta _{0}$ and $q_{**}$ is independent of $\beta _{\infty }$%
, it is convient to choose $\beta _{0}=\beta _{\infty }=0$, so that Theorems 
\ref{teo2} and \ref{teo4} yield to exponents $q_{1},q_{2}$ such that 
\[
1<q_{1}<q^{*}=2\frac{N-a+1}{N-4},\quad q_{2}>q_{**}=2\frac{2N-a-2}{2N-a-4}.
\]
If $a<4$, one has $q_{**}<q^{*}$ and therefore we get the compact
embedding 
\begin{equation}
H_{V,\mathrm{r}}^{1}\hookrightarrow L_{K}^{q}\qquad \textrm{for}\quad 
2\frac{2N-a-2}{2N-a-4}<q<2\frac{N-a+1}{N-4}.  \label{ES(sww): p}
\end{equation}
If $a=4$, then $q_{**}=q^{*}=2\left( N-3\right) /\left( N-4\right) $ and we
have the compact embedding 
\[
H_{V,\mathrm{r}}^{1}\hookrightarrow L_{K}^{q_{1}}+L_{K}^{q_{2}}\qquad \textrm{%
for}\quad 1<q_{1}<2\frac{N-3}{N-4}<q_{2}.
\]
Since $V$ and $K$ are power potentials, one can also apply the results of 
\cite{demarque-miya}, finding two exponents $s_{*}$ and $s^{*}$
such that the embedding $H_{V,\mathrm{r}}^{1}\hookrightarrow
L_{K}^{q}$ is compact if $s_{*}<q<s^{*}$. These exponents are exactly $q_{*}$ and $q^{*}$
of (\ref{qqq}) respectively, so that one obtains (\ref{ES(sww): p}) again,
provided that $a<4$. If $a=4$, instead, one gets $s_{*}=s^{*}$ and no result is avaliable in 
\cite{demarque-miya}.
\end{exa}

\section{Proof of Theorem \ref{teo1}}\label{thm1}
This section is devoted to the proof of Theorem \ref{teo1}, so let $N\geq 5$, assume \textsc{(\textbf{V})} and \textsc{(\textbf{K})} and take $q_1,q_2>1$. 
We begin with some preliminary results.

\begin{lem} \label{lemma1}
Take $R>r>0$ and $1<q<\infty$. Then there exist two constants $\tilde{C}=\tilde{C}(N,r,R,q,s)>0$ and $l=l(q,s)>0$ such that $\forall u\in H^2_{V,r}$ one has
\begin{equation} \label{integr3.3}
\int_{B_R\setminus B_r} K(|x|)|u|^q \, dx \leq \tilde{C} \left \| K(|\cdot|) \right \|_{L^s(B_R\setminus B_r)} \left \| u \right \|^{q-2l} \left( \int_{B_R\setminus B_r} |u|^2 \, dx \right)^l
\end{equation}
Furthermore, if $s>\frac{2N}{N+4}$ in assumption \textsc{(\textbf{K})}, then there exists $\tilde{C}_1=\tilde{C}_1(N,r,R,q,s)>0$ such that $\forall u\in H^2_{V,r}$ and $\forall h\in H^2_{V}$ we have
$$
\dfrac{\displaystyle\int_{B_R\setminus B_r} K(|x|)|u|^{q-1}|h| \, dx}{\tilde{C}_1 \left \| K(|\cdot|) \right \|_{L^s(B_R\setminus B_r)}} \leq \begin{cases}
\left( \displaystyle\int_{B_R\setminus B_r} |u|^2 \, dx \right)^{\frac{q-1}{2}} \left \| h \right \| & if \ q\leq \tilde{q} \\
\left( \displaystyle\int_{B_R\setminus B_r} |u|^2 \, dx \right)^{\frac{\tilde{q}-1}{2}} \left \| u \right \|^{q-\tilde{q}} \left \| h \right \| & if \ q>\tilde{q} \\
\end{cases}
$$
\noindent where $\tilde{q}:=2\left(1+\frac{2}{N}-\frac{1}{s} \right)$.
\end{lem}

\proof
Take $u\in H^2_{V,r}$ and fix $t\in (1,s)$ such that $t'q>2$ (where $t'=t/(t-1)$). Then, by H\"older inequality and \eqref{stimauso}, we get
\begin{equation*}
\begin{split}
\int_{B_R\setminus B_r}  K(|x|)|u|^q \, dx 
& \ \ \leq \left( \int_{B_R\setminus B_r} K(|x|)^t \, dx\right)^\frac{1}{t} \left( \int_{B_R\setminus B_r} |u|^{t'q} \, dx\right)^\frac{1}{t'} \\
& \ \ \leq |B_R\setminus B_r|^{\frac{1}{t}-\frac{1}{s}} \left\| K(|\cdot|) \right\|_{L^s(B_R\setminus B_r)} \left( \int_{B_R\setminus B_r} |u|^{t'q-2}|u|^2 \, dx\right)^\frac{1}{t'} \\
& \ \ \leq |B_R\setminus B_r|^{\frac{1}{t}-\frac{1}{s}} \left\| K(|\cdot|) \right\|_{L^s(B_R\setminus B_r)} \left( \dfrac{C_N \left\| u \right\|}{r^\frac{N-4}{2}} \right)^{q-\frac{2}{t'}} \left( \int_{B_R\setminus B_r} |u|^2 \, dx\right)^\frac{1}{t'}.
\end{split}
\end{equation*}
This proves \eqref{integr3.3}. 
To prove the second statement, take $ u\in H^2_{V,r}$ and $h\in H^2_{V}$. Let $\sigma:=\frac{2N}{N+4}$ be the H\"older conjugate exponent of $2^{**}$. Notice that $\frac{s}{\sigma}>1$. From H\"older inequality and \eqref{immsobo} we deduce
\begin{equation*}
\begin{split}
\int_{B_R\setminus B_r}  K(|x|)|u|^{q-1}|h| \, dx 
& \ \ \leq \left( \int_{B_R\setminus B_r}  K(|x|)^\sigma|u|^{(q-1)\sigma} \, dx\right)^\frac{1}{\sigma} \left( \int_{B_R\setminus B_r} |h|^{2^{**}} \, dx\right)^\frac{1}{2^{**}} \\
& \ \ \leq \left( \left( \int_{B_R\setminus B_r}  K(|x|)^s \, dx \right)^\frac{1}{\sigma} \left( \int_{B_R\setminus B_r} |u|^{(q-1)\sigma \left(\frac{s}{\sigma}\right)'} \, dx \right)^\frac{1}{\left(\frac{s}{\sigma}\right)'} \right)^\frac{1}{\sigma} S_N \left\| h \right\| \\
& \ \ \leq S_N \left\| K(|\cdot|) \right\|_{L^s(B_R\setminus B_r)} \left\| h \right\| \left( \int_{B_R\setminus B_r} |u|^{2\frac{q-1}{\tilde{q}-1}} \, dx\right)^\frac{\tilde{q}-1}{2},
\end{split}
\end{equation*}
thanks to the fact that $\sigma\left(\frac{s}{\sigma}\right)'=\frac{2Ns}{(N+4)s-2N}=\frac{2}{\tilde{q}-1}$. If $q\leq \tilde{q}$ we have
\begin{equation*}
\begin{split}
\int_{B_R\setminus B_r}  K(|x|)|u|^{q-1}|h| \, dx 
& \ \ \leq S_N \left\| K(|\cdot|) \right\|_{L^s(B_R\setminus B_r)} \left\| h \right\| \left( |B_R\setminus B_r|^{1-\frac{q-1}{\tilde{q}-1}} \left( \int_{B_R\setminus B_r} |u|^2 \, dx\right)^\frac{q-1}{\tilde{q}-1} \right)^\frac{\tilde{q}-1}{2} \\
& \ \ = S_N \left\| K(|\cdot|) \right\|_{L^s(B_R\setminus B_r)} \left\| h \right\| |B_R\setminus B_r|^\frac{\tilde{q}-q}{2} \left( \int_{B_R\setminus B_r} |u|^2 \, dx\right)^\frac{q-1}{2}.
\end{split}
\end{equation*}
On the other hand, if $q>\tilde{q}$, from \eqref{stimauso} we get
\begin{equation*}
\begin{split}
\int_{B_R\setminus B_r}  K(|x|)|u|^{q-1}|h| \, dx 
& \ \ \leq S_N \left\| K(|\cdot|) \right\|_{L^s(B_R\setminus B_r)} \left\| h \right\| \left( \int_{B_R\setminus B_r} |u|^{2\frac{q-1}{\tilde{q}-1}-2} |u|^2 \, dx \right)^\frac{\tilde{q}-1}{2} \\
& \ \ = S_N \left\| K(|\cdot|) \right\|_{L^s(B_R\setminus B_r)} \left\| h \right\| \left( \dfrac{C_N \left\| u \right\|}{r^\frac{N-4}{2}} \right)^{q-\tilde{q}} \left( \int_{B_R\setminus B_r} |u|^2 \, dx\right)^\frac{\tilde{q}-1}{2}.
\end{split}
\end{equation*}
Hence the thesis follows.
\endproof

We will also need the following result about the convergence in $L^{q_1}_K+L^{q_2}_K$. It is proved in \cite{bad4}. 

\begin{lem} \label{L4}
Let $\{ u_n \} \subseteq L^{q_1}_K+L^{q_2}_K$ be a sequence such that $\forall \epsilon>0$ there are $n_\epsilon>0$ and a sequence of measurable sets $E_{\epsilon,n}\subseteq\mathbb{R}^N$ satisfying
$$\forall n>n_\epsilon, \quad\int_{E_{\epsilon,n}} K(|x|)|u_n|^{q_1} \, dx + \int_{E_{\epsilon,n}^c} K(|x|)|u_n|^{q_2} \, dx<\epsilon.$$
Then $u_n\rightarrow 0$ in $L^{q_1}_K+L^{q_2}_K$.
\end{lem}

We can now prove Theorem \ref{teo1}. The arguments are similar to those of \cite{bad1}, so we will skip the details.

\proof[Proof of Theorem \ref{teo1}]
\begin{enumerate}
\item We can choose $R_1<R_2$ in hypothesis $\left( \mathcal{S}_{q_{1},q_{2}}^{\prime }\right) $. If $u\in H^2_{V,r}$, $u\neq 0$, we get
\begin{equation} \label{S0teo1}
\int_{B_{R_1}} K(|x|)|u|^{q_1} \, dx =\Vert u\Vert^{q_1}\int_{B_{R_1}} K(|x|)\dfrac{|u|^{q_1}}{\Vert u\Vert^{q_1}} \, dx 
\leq \Vert u\Vert^{q_1} \mathcal{S}_0(q_1,R_1)
\end{equation}
and similarly 
\begin{equation} \label{Sinfiteo1}
\int_{B_{R_2}^c} K(|x|)|u|^{q_2} \, dx \leq \Vert u\Vert^{q_2} \mathcal{S}_\infty(q_2,R_2).
\end{equation}
Furthemore, from Lemma \ref{lemma1} and the continuous embedding $D^{2,2}_r(\mathbb{R}^N)\hookrightarrow L^2_{\text{loc}}(\mathbb{R}^N)$, we obtain that there is a constant $C_1>0$, independent from $u$, such that
\begin{equation}\label{annulus}
\int_{B_{R_2}\setminus B_{R_1}} K(|x|)|u|^{q_1} \, dx \leq C_1\Vert u\Vert^{q_1}.
\end{equation}
Hence $u\in L^{q_1}_K(B_{R_2})\cap L^{q_2}_K(B^c_{R_2})$, so that $u\in L^{q_1}_K+ L^{q_2}_K$. 
Moreover, by \eqref{S0teo1}-\eqref{annulus} and Lemma \ref{L4}, one obtains that 
$u_n\rightarrow 0$ in $ H^2_{V,r}$ implies $u_n\rightarrow 0$ in $L^{q_1}_K+L^{q_2}_K$.

\item Assume $\left( \mathcal{S}_{q_{1},q_{2}}^{\prime\prime }\right) $, fix $\epsilon>0$ and choose a sequence $u_n\rightharpoonup 0$ in $H^2_{V,r}$. From \eqref{S0teo1}, \eqref{Sinfiteo1}, Lemma \ref{lemma1} and the compactness of the embedding $D^{2,2}_r(\mathbb{R}^N)\hookrightarrow L^2_{\text{loc}}(\mathbb{R}^N)$, we obtain
$$\int_{B_{R_\epsilon}} K(|x|)|u|^{q_1} \, dx + \int_{B_{R_\epsilon^c}} K(|x|)|u|^{q_2} \, dx < \epsilon$$
for any $n$ large enough. By Lemma \ref{L4}, this implies that $u_n\rightarrow0$ in $L^{q_1}_K+L^{q_2}_K$, which gives the compactness of the embedding.

\end{enumerate}
\endproof

\section{Proofs of Theorems \ref{teo2} and \ref{teo3} }\label{thms23}

\noindent In this section we let $N\geq 5$ and prove Theorems \ref{teo2} and \ref{teo3}. The first step is the following lemma, which will be also useful in the proofs of Theorems \ref{teo4} and \ref{teo5}.

\begin{lem} \label{lemma2}
Let $\Omega\subseteq\mathbb{R}^N$ a nonempty measurable set such that $V(|x|)<+\infty$ almost everywhere on $\Omega$
and assume that
$$\Lambda:=\esssup_{x\in \Omega} \ \dfrac{K(|x|)}{|x|^{\alpha} V(|x|)^{\beta}}<+\infty 
\quad \text{for some }\ 0\leq\beta\leq1 \text{ and }\alpha\in\mathbb{R}.$$
Take $u\in H^2_V$ and assume that there exist $\nu\in\mathbb{R}$ and $m>0$ such that
$$|u(x)|\leq\dfrac{m}{|x|^\nu} \quad \text{for almost every }\ x\in \Omega.$$
Then $\forall h\in H^2_V$ and $\forall q>\max\lbrace1,2\beta\rbrace$, we have
\begin{equation*}
 \int_\Omega K(|x|)|u|^{q-1}|h| \, dx 
 \leq\begin{cases} \Lambda m^{q-1} S^{1-2\beta}_N \left( \displaystyle\int_\Omega |x|^{\frac{\alpha-\nu(q-1)}{N+4(1-2\beta)}2N} \, dx \right)^{\frac{N+4(1-2\beta)}{2N}}  \left \| h \right \| & \text{if} \ \ 0 \le \beta \le \frac{1}{2} \\
\Lambda m^{q-2\beta}  \left( \displaystyle\int_\Omega |x|^{\frac{\alpha-\nu(q-2\beta)}{1-\beta}} \, dx \right)^{1-\beta} \left \| u \right \|^{2\beta-1}  \left \| h \right \| & \text{if} \ \ \frac{1}{2} < \beta < 1 \\
\Lambda m^{q-2}  \left( \displaystyle\int_\Omega |x|^{2\alpha-2\nu(q-2)} V(|x|)|u|^2 \, dx \right)^{\frac{1}{2}} \left \| h \right \| & \text{if} \ \ \beta=1.
\end{cases}
\end{equation*}
\end{lem}

\proof
We consider several cases.\smallskip

\noindent$\bullet$ Case $ \beta=0$. One has
\begin{equation*}
\begin{split}
\dfrac{1}{\Lambda}\int_\Omega K(|x|)|u|^{q-1}|h| \, dx  & \leq \int_\Omega |x|^\alpha |u|^{q-1}|h| \, dx 
\leq \left(\int_\Omega (|x|^\alpha |u|^{q-1})^\frac{2N}{N+4} \, dx \right)^\frac{N+4}{2N} \left(\int_\Omega |h|^{2^{**}} \, dx \right)^\frac{1}{2^{**}} \\
& \leq S_Nm^{q-1}\left( \int_\Omega |x|^{\frac{\alpha-\nu(q-1)}{N+4}2N} \, dx \right)^{\frac{N+4}{2N}}  \left \| h \right \|.
\end{split}
\end{equation*}

\noindent$\bullet$ Case $0<\beta<\frac{1}{2}$. We have $\frac{1}{\beta}>1$ and $\frac{1-\beta}{1-2\beta}2^{**}>1$ with H\"older conjugate exponents $(\frac{1}{\beta})'=\frac{1}{1-\beta}$ and $\left(\frac{1-\beta}{1-2\beta}2^{**}\right)'=\frac{2N(1-\beta)}{N+4(1-2\beta)}$. Then we get
\begin{equation*}
\begin{split}
\dfrac{1}{\Lambda}  \int_\Omega K(|x|)|u|^{q-1}|h| \, dx & \leq \int_\Omega |x|^\alpha V(|x|)^\beta |u|^{q-1}|h| \, dx \\
& \leq \left(\int_\Omega \left(|x|^\alpha |u|^{q-1}|h|^{1-2\beta}\right)^\frac{1}{1-\beta} \, dx \right)^{1-\beta} \left(\int_\Omega V(|x|)|h|^2 \, dx \right)^\beta \\
& \leq \left[\left(\int_\Omega \left(|x|^\alpha |u|^{q-1}\right)^{\frac{2N}{N+4(1-2\beta)}} \, dx \right)^\frac{N+4(1-2\beta)}{2N(1-\beta)} \left(\int_\Omega |h|^{2^{**}} \, dx \right)^\frac{(1-2\beta)}{(1-\beta)2^{**}}\right]^{1-\beta} \left \| h \right \|^{2\beta}\\
& \leq m^{q-1} \left(\int_\Omega |x|^{\frac{\alpha-\nu(q-1)}{N+4(1-2\beta)}2N} \, dx \right)^\frac{N+4(1-2\beta)}{2N} S^{1-2\beta}_N \left \| h \right \|^{1-2\beta} \left \| h \right \|^{2\beta}\\
& = m^{q-1} S^{1-2\beta}_N \left(\int_\Omega |x|^{\frac{\alpha-\nu(q-1)}{N+4(1-2\beta)}2N} \, dx \right)^\frac{N+4(1-2\beta)}{2N} \left \| h \right \|.
\end{split}
\end{equation*}

\noindent$\bullet$ Case $\beta=\frac{1}{2}$. We have
\begin{equation*}
\begin{split}
\dfrac{1}{\Lambda}\int_\Omega K(|x|)|u|^{q-1}|h| \, dx  & \leq \int_\Omega |x|^\alpha |u|^{q-1}V(|x|)^\frac{1}{2}|h| \, dx 
\leq \left( \int_\Omega |x|^{2\alpha} |u|^{2(q-1)} \, dx \right)^\frac{1}{2} \left(\int_\Omega V(|x|)|h|^2 \, dx \right)^\frac{1}{2} \\
& \leq m^{q-1}\left( \int_\Omega |x|^{2\alpha-2\nu(q-1)} \, dx \right)^{\frac{1}{2}}  \left \| h \right \|.
\end{split}
\end{equation*}

\noindent$\bullet$ Case $\ \frac{1}{2}<\beta<1$. We have $\frac{1}{2\beta-1}>1$ with H\"older conjugate exponent $\left(\frac{1}{2\beta-1}\right)'=\frac{1}{2(1-\beta)}$. Hence
\begin{equation*}
\begin{split}
\dfrac{1}{\Lambda}  \int_\Omega K(|x|)|u|^{q-1}|h| \, dx & \leq \int_\Omega |x|^\alpha V(|x|)^\beta |u|^{q-1}|h| \, dx \\
& \leq \left(\int_\Omega |x|^{2\alpha} V(|x|)^{2\beta-1} |u|^{2(q-1)} \, dx \right)^\frac{1}{2} \left(\int_\Omega V(|x|)|h|^2 \, dx \right)^\frac{1}{2} \\
& \leq \left(\int_\Omega |x|^{2\alpha} |u|^{2(q-2\beta)} V(|x|)^{2\beta-1} |u|^{2(2\beta-1)} \, dx \right)^\frac{1}{2} \left \| h \right \|\\
& \leq \left[\left(\int_\Omega |x|^\frac{\alpha}{1-\beta} |u|^\frac{q-2\beta}{1-\beta} \, dx \right)^{2(1-\beta)} \left(\int_\Omega V(|x|)|u|^2 \, dx \right)^{2\beta-1}\right]^\frac{1}{2} \left \| h \right \|\\
& \leq m^{q-2\beta} \left[ \left(\int_\Omega |x|^{\frac{\alpha}{1-\beta}-\nu\frac{q-2\beta}{1-\beta}} \, dx \right)^{2(1-\beta)} \left(\int_\Omega V(|x|)|u|^2 \, dx \right)^{2\beta-1} \right]^\frac{1}{2}  \left \| h \right \|\\
& \leq m^{q-2\beta} \left(\int_\Omega |x|^\frac{\alpha-\nu(q-2\beta)}{1-\beta} \, dx \right)^{1-\beta} \left \| u \right \|^{2\beta-1} \left \| h \right \|.
\end{split}
\end{equation*}

\noindent$\bullet$ Case $ \beta=1$. As $q>\max\{1,2\beta\}$, in this case we have $q>2$. Hence
\begin{equation*}
\begin{split}
\dfrac{1}{\Lambda}  \int_\Omega K(|x|)|u|^{q-1}|h| \, dx & \leq \int_\Omega |x|^\alpha V(|x|) |u|^{q-1}|h| \, dx \\
& \leq \left(\int_\Omega |x|^{2\alpha} V(|x|) |u|^{2(q-1)} \, dx \right)^\frac{1}{2} \left(\int_\Omega V(|x|)|h|^2 \, dx \right)^\frac{1}{2} \\
& \leq \left(\int_\Omega |x|^{2\alpha} |u|^{2(q-2)} V(|x|) |u|^{2} \, dx \right)^\frac{1}{2} \left \| h \right \|\\
& \leq m^{q-2} \left(\int_\Omega |x|^{2\alpha-2\nu(q-2)} V(|x|) |u|^{2} \, dx \right)^\frac{1}{2} \left \| h \right \|.
\end{split}
\end{equation*}
\endproof

We can now give the proofs of Theorems $\ref{teo2}$ and \ref{teo3}.

\proof[Proof of Theorem $\ref{teo2}$]
Let $u\in H^2_{V,r}$ and $h\in H^2_V$ such that $\left\| u \right\|=\left\| h \right\|=1$. Let $0<R \leq R_1$. 
Thanks to \eqref{stimauso} and
$$\esssup_{x\in B_R} \ \dfrac{K(|x|)}{|x|^{\alpha_0} V(|x|)^{\beta_0}} \ \leq \ \esssup_{r\in(0,R_1)} \ \dfrac{K(r)}{r^{\alpha_0} V(r)^{\beta_0}}<+\infty,$$
we can apply Lemma $ \ref{lemma2}$ with $\Omega=B_R$, $\alpha=\alpha_0$, $\beta=\beta_0$, $m=C_N$ and $\nu=\frac{N-4}{2}$. 
In the following $C$ is any positive constant independent from $u$, $h$ and $R$.
If $0\leq\beta_0\leq\frac{1}{2}$ we get
\begin{equation*}
\begin{split}
\int_{B_R} K(|x|)|u|^{q_1-1}|h| \, dx  & \leq C\left( \int_{B_R} |x|^{\frac{\alpha_0-\frac{N-4}{2}(q_1-1)}{N+4(1-2\beta_0)}2N} \, dx \right)^{\frac{N+4(1-2\beta_0)}{2N}}\\
& \leq C\left( \int_0^R r^{\frac{2\alpha_0-(N-4)(q_1-1)}{N+4(1-2\beta_0)}N+N-1} \, dr \right)^{\frac{N+4(1-2\beta_0)}{2N}}\\
& = C\left( R^{\frac{2\alpha_0-8\beta_0+2N-(N-4)q_1}{N+4(1-2\beta_0)}N} \right)^{\frac{N+4(1-2\beta_0)}{2N}},
\end{split}
\end{equation*}
because
\[
2\alpha_0-8\beta_0+2N-(N-4)q_1 \\ =(N-4)\left( 2\dfrac{\alpha_0-4\beta_0+N}{N-4}-q_1\right)
 =(N-4)\left(q^*(\alpha_0,\beta_0)-q_1\right)>0.
\]
If $\frac{1}{2}<\beta_0<1$ we get
\begin{equation*}
\begin{split}
\int_{B_R} K(|x|)|u|^{q_1-1}|h| \, dx  & \leq C\left( \int_{B_R} |x|^{\frac{\alpha_0-\frac{N-4}{2}(q_1-2\beta_0)}{1-\beta_0}} \, dx \right)^{1-\beta_0}
\leq C\left( \int_0^R r^{\frac{2\alpha_0-(N-4)(q_1-2\beta_0)}{2(1-\beta_0)}+N-1} \, dr \right)^{1-\beta_0}\\
& = C\left( R^{\frac{2\alpha_0-(N-4)(q_1-2\beta_0)}{2(1-\beta_0)}+N} \right)^{1-\beta_0},
\end{split}
\end{equation*}
because
\begin{gather*}
{\dfrac{2\alpha_0-(N-4)(q_1-2\beta_0)}{2(1-\beta_0)}+N} 
=\dfrac{(N-4)}{2(1-\beta_0)}\left( 2\dfrac{\alpha_0-4\beta_0+N}{N-4}-q_1\right)=\dfrac{(N-4)}{2(1-\beta_0)}\left(q^*(\alpha_0,\beta_0)-q_1\right)>0.
\end{gather*}
If $\beta_0=1$ we get
\begin{equation*}
\begin{split}
\int_{B_R} K(|x|)|u|^{q_1-1}|h| \, dx  & \leq C\left( \int_{B_R} |x|^{2\alpha_0-(N-4)(q_1-2)}V(|x|)|u|^2 \, dx \right)^\frac{1}{2}\\
& \leq C\left( R^{2\alpha_0-(N-4)(q_1-2)}\int_{B_R} V(|x|)|u|^2 \, dx \right)^\frac{1}{2}
 \leq C\ R^\frac{2\alpha_0-(N-4)(q_1-2)}{2},\\
\end{split}
\end{equation*}
because
\begin{gather*}
2\alpha_0-8+2N-(N-4)q_1 
=(N-4)\left( 2\dfrac{\alpha_0-4+N}{N-4}-q_1\right)=(N-4)\left(q^*(\alpha_0,1)-q_1\right)>0.
\end{gather*}
Hence, in any case, we get $\mathcal{R}_0(q_1,R)\leq CR^\delta$ for some $\delta=\delta(N,\alpha_0,\beta_0,q_1)>0$, which gives the result.
\endproof

\proof[Proof of Theorem $\ref{teo3}$]
Let $u\in H^2_{V,r}$ and $h\in H^2_V$ such that $\left\| u \right\|=\left\| h \right\|=1$. Let $R \geq R_2$. By \eqref{stimauso} and
 $$\esssup_{x\in B_R^C} \ \dfrac{K(|x|)}{|x|^{\alpha_\infty} V(|x|)^{\beta_\infty}} \ \leq \ \esssup_{r>R_2} \ \dfrac{K(r)}{r^{\alpha_\infty} V(r)^{\beta_\infty}}<+\infty,$$
we can apply lemma $ \ref{lemma2}$ with $\Omega=B_R^c$, $\alpha=\alpha_\infty$, $\beta=\beta_\infty$, $m=C_N$ and $\nu=\frac{N-4}{2}$. Hereafter $C$ denotes any positive constant independent form $u$, $h$ and $R$.
If $0\leq\beta_\infty\leq\frac{1}{2}$ we get
\begin{equation*}
\begin{split}
\int_{B_R^c} K(|x|)|u|^{q_2-1}|h| \, dx  & \leq C\left( \int_{B_R^c} |x|^{\frac{\alpha_\infty-\frac{N-4}{2}(q_2-1)}{N+4(1-2\beta_\infty)}2N} \, dx \right)^{\frac{N+4(1-2\beta_\infty)}{2N}}\\
& \leq C\left( \int_R^{+\infty} r^{\frac{2\alpha_\infty-(N-4)(q_2-1)}{N+4(1-2\beta_\infty)}N+N-1} \, dr \right)^{\frac{N+4(1-2\beta_\infty)}{2N}}\\
& = C\left( R^{\frac{2\alpha_\infty-8\beta_\infty+2N-(N-4)q_2}{N+4(1-2\beta_\infty)}N} \right)^{\frac{N+4(1-2\beta_\infty)}{2N}},
\end{split}
\end{equation*}
because
\begin{gather*}
2\alpha_\infty-8\beta_\infty+2N-(N-4)q_2 
=(N-4)\left( 2\dfrac{\alpha_\infty-4\beta_\infty+N}{N-4}-q_2\right)
=(N-4)\left(q^*(\alpha_\infty,\beta_\infty)-q_2\right)<0.
\end{gather*}
If $\frac{1}{2}<\beta_\infty<1$ we have
\begin{equation*}
\begin{split}
\int_{B_R^c} K(|x|)|u|^{q_2-1}|h| \, dx  & \leq C\left( \int_{B_R^c} |x|^{\frac{\alpha_\infty-\frac{N-4}{2}(q_2-2\beta_\infty)}{1-\beta_\infty}} \, dx \right)^{1-\beta_\infty}\\
& \leq C\left( \int_R^{+\infty} r^{\frac{2\alpha_\infty-(N-4)(q_2-2\beta_\infty)}{2(1-\beta_\infty)}+N-1} \, dr \right)^{1-\beta_\infty}\\
& = C\left( R^{\frac{2\alpha_\infty-(N-4)(q_2-2\beta_\infty)}{2(1-\beta_\infty)}+N} \right)^{1-\beta_\infty},
\end{split}
\end{equation*}
because
\begin{equation*}
\begin{split}
{\dfrac{2\alpha_\infty-(N-4)(q_2-2\beta_\infty)}{2(1-\beta_\infty)}+N} 
& = \dfrac{(N-4)}{2(1-\beta_\infty)}\left( 2\dfrac{\alpha_\infty-4\beta_\infty+N}{N-4}-q_2\right)\\
& = \dfrac{(N-4)}{2(1-\beta_\infty)}\left(q^*(\alpha_\infty,\beta_\infty)-q_2\right)<0.
\end{split}
\end{equation*}
If $\beta_\infty=1$ we get
\begin{equation*}
\begin{split}
\int_{B_R^c} K(|x|)|u|^{q_2-1}|h| \, dx  & \leq C\left( \int_{B_R^c} |x|^{2\alpha_\infty-(N-4)(q_2-2)}V(|x|)|u|^2 \, dx \right)^\frac{1}{2}\\
& \leq C\left( R^{2\alpha_\infty-(N-4)(q_2-2)}\int_{B_R^C} V(|x|)|u|^2 \, dx \right)^\frac{1}{2}
 \leq C\ R^\frac{2\alpha_\infty-(N-4)(q_2-2)}{2},\\
\end{split}
\end{equation*}
because
\begin{gather*}
2\alpha_\infty-8+2N-(N-4)q_2 
=(N-4)\left( 2\dfrac{\alpha_\infty-4+N}{N-4}-q_2\right)=(N-4)\left(q^*(\alpha_\infty,1)-q_2\right)<0.
\end{gather*}
So, in any case, we get $\mathcal{R}_\infty(q_2,R)\leq CR^\delta$ for some $\delta=\delta(N,\alpha_0,\beta_0,q_1)<0$. Hence the thesis follows.
\endproof

\section{Proofs of Theorems \ref{teo4} and \ref{teo5}}\label{thms45}

Let $N\geq 5$. To prove Theorems \ref{teo4} and \ref{teo5} we need some preliminary results about pointwise estimates of radial Sobolev functions.

For any open interval $\mathcal{I}\subset\mathbb{R}$, we will consider the space
$$W^{2,1}(\mathcal{I}):=\left\lbrace u\in L^1(\mathcal{I}): D^\alpha u\in L^1(\mathcal{I}), \ \forall |\alpha|\leq2 \right\rbrace.$$
The proof of the following lemma can be easily derived from the arguments of \cite[Appendix]{bad2}, so we skip it.

\begin{lem} \label{lemminteg}
Let $u \in D_{r}^{2,2}(\mathbb{R}^N)$ and define $\tilde{u}:\mathbb{R}_+\rightarrow\mathbb{R}$ such that $u(x)=\tilde{u}(|x|)$ for almost every $x\in \mathbb{R}^N$. Then $\tilde{u}\in W^{2,1}(\mathcal{I})$ for every open bounded interval $\mathcal{I}\subset\mathbb{R}_+$ such that $\inf\mathcal{I}>0$.
\end{lem}

\begin{prop} \label{radfuori}
Assume that there exists $R_2>0$ such that $V(r)<+\infty$ for almost every $r>R_2$ and
$$\lambda_\infty:=\essinf_{r>R_2} \ r^{\gamma_\infty}V(r)>0 \quad \text{for some } \gamma_\infty\leq \dfrac{14}{3}.$$
Then $\forall u\in H^2_{V,r}$ we have
\begin{equation} \label{eqfuori}
|u(x)|\leq c_\infty \lambda_\infty^{-\frac{1}{4}}\dfrac{\left \| u \right \|}{|x|^\frac{2(N-2)-\gamma_\infty}{4}} 
\quad \text{almost everywhere in } B^c_{R_2},
\end{equation}
\noindent where $c_\infty=\dfrac{1}{\sqrt{\sigma_N}}\left(\dfrac{8}{N(2(N-2)-\gamma_\infty)}\right)^\frac{1}{4}$.
\end{prop}

Notice that $2(N-2)-\gamma_\infty>0$ in \eqref{eqfuori}.

\proof
\noindent Let $u\in H^2_{V,r}$. Define $\tilde{u}:\mathbb{R}_+\rightarrow\mathbb{R}$ as the continuous function such that $u(x)=\tilde{u}(|x|)$ for almost every $x\in \mathbb{R}^N$. Define
$$v(r):=r^{\frac{N}{2}-1-\frac{\gamma_\infty}{4}}\tilde{u}(r)^2 \quad\text{for every} \ r>0.$$

\noindent If $\lambda:=\liminf_{r\rightarrow+\infty}v(r)>0$, then for $r$ large enough we get
$$r^{N-1-\gamma_\infty}\tilde{u}(r)^2\geq \dfrac{\lambda}{2r^{\frac{3}{4}\gamma_\infty-\frac{N}{2}}}$$
and from this we get the following contradiction:

\begin{equation*}
\begin{split}
\int_{B^c_{R_2}}V(|x|)u^2 \, dx 
 \geq \lambda_\infty \int_{B^c_{R_2}}\dfrac{u^2}{|x|^{\gamma_\infty}} \, dx =\lambda_\infty \sigma_N \int_{R_2}^{+\infty}\dfrac{\tilde{u}(r)^2}{r^{\gamma_\infty}} r^{N-1} \, dr 
 \geq \lambda_\infty \sigma_N \int_{R_2}^{+\infty}\dfrac{\lambda}{2r^{\frac{3}{4}\gamma_\infty-\frac{N}{2}}} \, dr=+\infty
\end{split}
\end{equation*}
where the last integral diverges because $N\geq5$ and $\gamma_\infty\leq14/3$. Hence, it must be $\lambda=0$ and therefore there is a sequence $r_n\rightarrow+\infty$ such that $v(r_n)\rightarrow\nobreak0$.
From lemma \ref{lemminteg} we get $v\in W^{2,1}((r,r_n))$ for all $R_2<r<r_n<+\infty$, whence
$$v(r_n)-v(r)=\int_r^{r_n}v'(s)\, ds.$$
Furthermore, for all $s\in (r,r_n)$ we have
\begin{equation*}
\begin{split}
v'(s) &=\left( \dfrac{N}{2}-1-\dfrac{\gamma_\infty}{4} \right)s^{\frac{N}{2}-2-\frac{\gamma_\infty}{4}}\tilde{u}(s)^2+2s^{\frac{N}{2}-1-\frac{\gamma_\infty}{4}}\tilde{u}(s)\tilde{u}'(s) 
\geq 2s^{\frac{N}{2}-1-\frac{\gamma_\infty}{4}}\tilde{u}(s)\tilde{u}'(s) \\
& \geq -2s^{\frac{N}{2}-1-\frac{\gamma_\infty}{4}}|\tilde{u}(s)||\tilde{u}'(s)|.
\end{split}
\end{equation*}
The first inequality derives from $\frac{N}{2}-1-\frac{\gamma_\infty}{4}>\frac{N}{2}-\frac{13}{6}>0$. Then we get
$$v(r_n)-v(r)=\int_r^{r_n}v'(s)\, ds \geq -2\int_r^{r_n}s^{\frac{N}{2}-1-\frac{\gamma_\infty}{4}}|\tilde{u}(s)||\tilde{u}'(s)| \, ds.$$
Now from \eqref{stima2} we deduce that
\begin{equation*}
\begin{split}
v(r)-v(r_n) &\leq 2\int_r^{r_n}s^{\frac{N}{2}-1-\frac{\gamma_\infty}{4}}|\tilde{u}(s)||\tilde{u}'(s)|\, ds 
\leq \dfrac{2}{\sqrt{N\sigma_N}}\int_r^{r_n}s^{\frac{N}{2}-1-\frac{\gamma_\infty}{4}}|\tilde{u}(s)| \dfrac{\left \| \Delta u \right \|_2}{s^\frac{N-2}{2}} \, ds \\
&= \dfrac{2\left \| \Delta u \right \|_2}{\sqrt{N\sigma_N}}\int_r^{r_n}\dfrac{|\tilde{u}(s)|}{s^\frac{\gamma_\infty}{2}}s^\frac{N-1}{2} \dfrac{1}{s^{\frac{N}{2}-\frac{1}{2}-\frac{\gamma_\infty}{4}}} \, ds \\
& \leq \dfrac{2\left \| \Delta u \right \|_2}{\sqrt{N\sigma_N}} \left( \int_r^{r_n}\dfrac{\tilde{u}(s)^2}{s^{\gamma_\infty}}s^{N-1} \, ds \right)^{\frac{1}{2}} \left(\int_r^{r_n} \dfrac{1}{s^{N-1-\frac{\gamma_\infty}{2}}} \, ds \right)^\frac{1}{2}\\
& \leq \dfrac{2\left \| \Delta u \right \|_2}{\sqrt{N\sigma_N}} \left( \int_{R_2}^{+\infty}\dfrac{\tilde{u}(s)^2}{s^{\gamma_\infty}}s^{N-1} \, ds \right)^{\frac{1}{2}} \left(\int_r^{+\infty} \dfrac{1}{s^{2 \left( \frac{N}{2}-1-\frac{\gamma_\infty}{4} \right)+1}} \, ds \right)^\frac{1}{2}\\
& \leq \dfrac{2\left \| \Delta u \right \|_2}{\sqrt{N\sigma_N}} \left( \dfrac{1}{\lambda_\infty}\int_{R_2}^{+\infty}V(s)\tilde{u}(s)^2s^{N-1} \, ds \right)^{\frac{1}{2}} \left( \left[ \dfrac{s^{-2 \left( \frac{N}{2}-1-\frac{\gamma_\infty}{4} \right)}}{-2\left( \frac{N}{2}-1-\frac{\gamma_\infty}{4} \right) } \right]_r^{+\infty}\right)^\frac{1}{2}\\
& \leq \dfrac{2\left \| \Delta u \right \|_2}{\sqrt{N\sigma_N}} \left( \dfrac{1}{\lambda_\infty\sigma_N}\int_{\mathbb{R}^N}V(|x|)u^2 \, dx \right)^{\frac{1}{2}}  \dfrac{r^{- \left( \frac{N}{2}-1-\frac{\gamma_\infty}{4} \right)}}{\sqrt{2\left( \frac{N}{2}-1-\frac{\gamma_\infty}{4} \right)}}.
\end{split}
\end{equation*}
Since it is easy to see that $\left \| \Delta u \right \|_2\left \| u \right \|_V\leq \left \| u \right \|^2$, we obtain
$$v(r)-v(r_n)\leq  \dfrac{1}{\sigma_N}\sqrt{\dfrac{8}{\lambda_\infty N[2(N-2)-\gamma_\infty]}} \dfrac{\left \| u \right \|^2}{r^{ \left( \frac{N}{2}-1-\frac{\gamma_\infty}{4} \right)}}.$$
Finally, recalling the definition of $v(r)$ e the fact that $v(r_n)\rightarrow 0$, we conclude
$$|x|^{ \left( \frac{N}{2}-1-\frac{\gamma_\infty}{4} \right)} |u(x)|^2\leq  \frac{1}{\sigma_N}\sqrt{\dfrac{8}{\lambda_\infty N[2(N-2)-\gamma_\infty]}} \dfrac{\left \| u \right \|^2}{|x|^{ \left( \frac{N}{2}-1-\frac{\gamma_\infty}{4} \right)}}$$
and hence
$|u(x)|\leq c_\infty \lambda_\infty^{-\frac{1}{4}} \left \| u \right \|\,|x|^{-\frac{2(N-2)-\gamma_\infty}{4}}.$
\endproof

We now prove a second pointwise estimate.

\begin{prop} \label{raddentro}
Assume there exists $R>0$ sucht that $V(r)<+\infty$ almost everywhere on $(0,R)$ and
$$\lambda_0:=\essinf_{r\in(0,R)} \ r^{\gamma_0}V(r)>0 \quad \text{for some } \gamma_0\geq 4.$$
Then $\forall u\in H^2_{V,r}$ we have
\begin{equation} \label{eqdentro}
|u(x)| \leq c_0\left( \dfrac{1}{\sqrt{\lambda_0}}+\dfrac{R^\frac{\gamma_0-4}{2}}{\lambda_0} \right)^\frac{1}{2} \dfrac{\left \|  u \right \|}{|x|^{\frac{2N-4-\gamma_0}{2}}} \quad \text{almost everywhere in} \ B_{R},
\end{equation}
where $c_0=\sqrt{\dfrac{\max\{2/\sqrt{N},N-7/2\}}{\sigma_N}}$.
\end{prop}

\proof
Let $u\in H^2_{V,r}$ and define $\tilde{u}:\mathbb{R}_+\rightarrow\mathbb{R}$ as the continuous function such that
$u(x)=\tilde{u}(|x|)$ for almost every $x\in \mathbb{R}^N$. 
Define
$$v(r):=r^{\frac{2N-3-\gamma_0}{2}}\tilde{u}(r)^2 \quad\text{for all} \ r>0.$$
If $\lambda:=\liminf_{r\rightarrow+\infty}v(r)>0$, then for all $r$ large enough we have
$$r^{N-1-\gamma_0}\tilde{u}(r)^2\geq \dfrac{\lambda}{2r^{\frac{\gamma_0}{2}-\frac{1}{2}}},$$
from which we derive a contradiction as follows:
\begin{equation*}
\begin{split}
\int_{B^{R}}V(|x|)u^2 \, dx & \geq \lambda_\infty \int_{B_{R}}\dfrac{u^2}{|x|^{\gamma_0}} \, dx =\lambda_0 \sigma_N \int_0^R\dfrac{\tilde{u}(r)^2}{r^{\gamma_0}} r^{N-1} \, dr 
\geq \lambda_0 \sigma_N \int_0^R\dfrac{\lambda}{2r^{\frac{\gamma_0-1}{2}}} \, dr=+\infty
\end{split}
\end{equation*}
where the last integral diverges because $\gamma_0\geq4$.
This proves that $\lambda=0$ and thus implies that there exists a sequence $r_n\rightarrow0^+$ tale che $v(r_n)\rightarrow\nobreak0$.
From lemma \ref{lemminteg} we get $v\in W^{2,1}((r_n,r))$ for all $0<r_n<r<R$, whence
$$v(r)-v(r_n)=\int^r_{r_n}v'(s)\, ds.$$
Furthemore for all $s\in (r_n,r)$ we have

$$v'(s)=\left( \dfrac{2N-3-\gamma_0}{2} \right)s^{\frac{2N-5-\gamma_0}{2}}\tilde{u}(s)^2+2 s^{\frac{2N-3-\gamma_0}{2}}\tilde{u}(s)\tilde{u}'(s)
=\left( \dfrac{2N-3-\gamma_0}{2} \right)I(s)+2 J(s)$$
with obvious definitions of $I(s)$ and $J(s)$, on which we obtain the following estimates:
\begin{equation*}
\begin{split}
\int_{r_n}^r I(s)  ds & = \int_{r_n}^r s^{\frac{2N-5-\gamma_0}{2}}\tilde{u}(s)^2 \, ds = \int_{r_n}^r \dfrac{\tilde{u}(s)^2}{s^{\gamma_0}}s^{N-1}s^\frac{\gamma_0-3}{2} \, ds 
 \leq r^\frac{\gamma_0-3}{2}\int_{0}^R \dfrac{\tilde{u}(s)^2}{s^{\gamma_0}}s^{N-1} \, ds \\
& \leq r^\frac{\gamma_0-3}{2} \dfrac{1}{\lambda_0} \int_{0}^R V(s) \tilde{u}(s)^2 s^{N-1} \, ds 
 \leq r^\frac{\gamma_0-3}{2} \dfrac{\left \| u \right \|^2}{\lambda_0 \sigma_N} \leq R^\frac{\gamma_0-4}{2} \dfrac{\left \| u \right \|^2}{\lambda_0 \sigma_N} r^\frac{1}{2}
\end{split}
\end{equation*}
and, by \eqref{stima2},
\begin{equation*}
\begin{split}
\int_{r_n}^r J(s)  ds & = \int_{r_n}^r s^{\frac{2N-3-\gamma_0}{2}}\tilde{u}(s)\tilde{u}'(s) \, ds \leq \int_{r_n}^r s^{\frac{2N-3-\gamma_0}{2}} |\tilde{u}(s)||\tilde{u}'(s)| \, ds \\
& \leq \dfrac{\left \| \Delta u \right \|_2}{\sqrt{N\sigma_N}}\int_{r_n}^r s^{\frac{2N-3-\gamma_0}{2}}|\tilde{u}(s)|\dfrac{1}{s^\frac{N-2}{2}} \, ds 
 = \dfrac{\left \| \Delta u \right \|_2}{\sqrt{N\sigma_N}}\int_{r_n}^r \dfrac{|\tilde{u}(s)|}{s^\frac{\gamma_0}{2}}s^\frac{N-1}{2} \, ds \\
& \leq  \dfrac{\left \| \Delta u \right \|_2}{\sqrt{N\sigma_N}} \left(\int_{r_n}^r \dfrac{\tilde{u}(s)^2}{s^{\gamma_0}}s^{N-1} \, ds \right)^\frac{1}{2} \left( \int_{r_n}^r \, ds \right)^\frac{1}{2} \\
& \leq \dfrac{\left \| \Delta u \right \|_2}{\sqrt{N\sigma_N}} \left( \dfrac{1}{\lambda_0} \int_{0}^R V(s) \tilde{u}(s)^2 s^{N-1} \, ds \right)^\frac{1}{2} \left( \int_{0}^r \, ds \right)^\frac{1}{2} \\
& \leq \dfrac{\left \| \Delta u \right \|_2}{\sigma_N\sqrt{N}}\dfrac{\left \|  u \right \|_V}{\sqrt{\lambda_0}} r^\frac{1}{2} \leq \dfrac{\left \|  u \right \|^2}{\sigma_N\sqrt{N\lambda_0}} r^\frac{1}{2}.
\end{split}
\end{equation*}
Now, if $4\leq \gamma_0 \leq 2N-3$, we get
\begin{equation*}
\begin{split}
v(r)-v(r_n) &= \int_{r_n}^r v'(s) \, ds \leq\left(\dfrac{2N-3-\gamma_0}{2}\right)\int_{r_n}^r I(s) ds +2\int_{r_n}^r J(s) ds \\
& \leq \left(\dfrac{2N-3-\gamma_0}{2}\right) R^\frac{\gamma_0-4}{2} \dfrac{\left \| u \right \|}{\lambda_0 \sigma_N} r^\frac{1}{2}+\dfrac{2\left \|  u \right \|^2}{\sigma_N\sqrt{N\lambda_0}} r^\frac{1}{2} \\
& \leq \left(N-\dfrac{7}{2}\right) R^\frac{\gamma_0-4}{2} \dfrac{\left \| u \right \|^2}{\lambda_0 \sigma_N} r^\frac{1}{2}+\dfrac{2\left \|  u \right \|^2}{\sigma_N\sqrt{N\lambda_0}} r^\frac{1}{2} 
\end{split}
\end{equation*}
On the other hand, if $\gamma_0 \geq 2N-3$, we get
$$v'(s)=\left( \dfrac{2N-3-\gamma_0}{2} \right)I(s)+2J(s) \leq 2J(s)$$
and thus
$$v(r)-v(r_n) = \int_{r_n}^r v'(s) \, ds \leq2\int_{r_n}^r J(s) ds \leq \dfrac{2\left \|  u \right \|^2}{\sigma_N\sqrt{N\lambda_0}} r^\frac{1}{2}.$$
So, in any case, we have
$$v(r)-v(r_n)\leq \dfrac{1}{\sigma_N}\left[ \dfrac{2}{\sqrt{N\lambda_0}}+\left( N-\dfrac{7}{2}\right)\dfrac{R^\frac{\gamma_0-4}{2}}{\lambda_0} \right] \left \|  u \right \|^2 r^\frac{1}{2}.$$
Hence, recalling the definition of $v(r)$ and the fact that $v(r_n)\rightarrow0$, we deduce
$$|x|^{\frac{2N-3-\gamma_0}{2}} |u(x)|^2 \leq \dfrac{1}{\sigma_N}\left[ \dfrac{2}{\sqrt{N\lambda_0}}+\left( N-\dfrac{7}{2}\right)\dfrac{R^\frac{\gamma_0-4}{2}}{\lambda_0} \right] \left \|  u \right \|^2 |x|^\frac{1}{2},$$
which gives \eqref{eqdentro}.
\endproof

We will also need the following lemma.

\begin{lem} \label{lemma5}
Assume that there exists $R>0$ be such that $V(r)<+\infty$ almost everywhere on $(0,R)$ and
\begin{equation} \label{lamgran}
\Lambda_{\alpha,\beta}(R):=\esssup_{r\in (0,R)} \ \dfrac{K(r)}{r^\alpha V(r)^\beta}<+\infty 
\quad\text{for some } \dfrac{1}{2}\leq \beta \leq 1 \text{ and } \alpha\in \mathbb{R}
\end{equation}
and
$$\lambda(R):=\essinf_{r\in(0,R)} \ r^{\gamma_0}V(r)>0 \quad\text{for some }  \gamma_0>4.$$
Assume also that $\exists q>2\beta$ such that
$(2N-4-\gamma_0)q<4\alpha+4N-2(\gamma_0+4)\beta.$
Then $\forall u\in H^2_{V,r}$ and $\forall h\in H^2_{V}$ we have
$$\int_{B_R} K(|x|)|u|^{q-1}|h| \, dx \leq c_0^{q-2\beta} a(R) R^\frac{4\alpha+4N-2(\gamma_0+4)\beta-(2N-4-\gamma_0)q}{4}\left \| u \right \|^{q-1} \left \| h \right \|,$$
where 
$a(R):=\Lambda_{\alpha,\beta}(R)\left( \dfrac{1}{\sqrt{\lambda(R)}}+\dfrac{R^\frac{\gamma_0-4}{2}}{\lambda(R)} \right)
^{(q-2\beta)/2}$ 
and $c_0$ is given in Proposition \ref{raddentro}.
\end{lem}

\proof
Take $u\in  H^2_{V,r}$ 
and $h\in H^2_{V}$. Thanks to assumption \eqref{lamgran} and Proposition \ref{raddentro}, we can apply Lemma \ref{lemma2} with $\Omega=B_R$, $\Lambda=\Lambda_{\alpha,\beta}(R)$, $\nu=\frac{2N-4-\gamma_0}{4}$ and
$$m=c_0\left( \dfrac{1}{\sqrt{\lambda(R)}}+\dfrac{R^\frac{\gamma_0-4}{2}}{\lambda(R)} \right)^\frac{1}{2} \left \| u \right \|.$$
If $\frac{1}{2}\leq\beta<1$ we have
\begin{equation*}
\begin{split}
\int_{B_R}  K(|x|)|u|^{q-1}|h| \, dx 
& \leq \Lambda m^{q-2\beta} \left( \int_\Omega |x|^{\frac{\alpha-\nu(q-2\beta)}{1-\beta}} \, dx \right)^{1-\beta} \left \| u \right \|^{2\beta-1}  \left \| h \right \| \\
& =c_0^{q-2\beta} a(R) \left( \int_{B_R} |x|^{\frac{4\alpha-(2N-4-\gamma_0)(q-2\beta)}{4(1-\beta)}} \, dx \right)^{1-\beta} \left \| u \right \|^{2\beta-1} \left \| u \right \|^{q-2\beta} \left \| h \right \| \\
& \leq c_0^{q-2\beta} a(R) \left( R^{\frac{4\alpha-(2N-4-\gamma_0)(q-2\beta)}{4(1-\beta)}+N} \right)^{1-\beta} \left \| u \right \|^{q-1} \left \| h \right \|,
\end{split}
\end{equation*}
because
$$\dfrac{4\alpha-(2N-4-\gamma_0)(q-2\beta)}{4(1-\beta)}+N=\dfrac{4\alpha+4N-2(\gamma_0+4)\beta-(2N-4-\gamma_0)q}{4(1-\beta)} >0$$
If $\beta=1$ we get
\begin{equation*}
\begin{split}
\int_{B_R}  K(|x|)|u|^{q-1}|h| \, dx 
& \leq \Lambda m^{q-2} \left( \int_\Omega |x|^{2\alpha-2\nu(q-2)}V(|x|)|u|^2 \, dx \right)^\frac{1}{2} \left \| h \right \| \\
& =c_0^{q-2\beta} a(R) \left( \int_{B_R} |x|^{\frac{4\alpha-(2N-4-\gamma_0)(q-2)}{2}} V(|x|)|u|^2 \, dx \right)^\frac{1}{2} \left \| u \right \|^{q-2} \left \| h \right \| \\
& \leq c_0^{q-2\beta} a(R) \left( R^{\frac{4\alpha-(2N-4-\gamma_0)(q-2)}{2}}\int_{B_R}  V(|x|)|u|^2 \, dx \right)^\frac{1}{2} \left \| u \right \|^{q-2} \left \| h \right \| \\
& \leq c_0^{q-2\beta} a(R) R^{\frac{4\alpha-(2N-4-\gamma_0)(q-2)}{4}}\left \| u \right \|^{q-1} \left \| h \right \|,
\end{split}
\end{equation*}
because $4\alpha-(2N-4-\gamma_0)(q-2)=4\alpha+4N-2(\gamma_0+2)-(2N-4-\gamma_0)q>0.$
\endproof

We can now give the proofs of Theorems \ref{teo4} and \ref{teo5}.
For convenience, define three functions $\alpha_1=\alpha_1(\beta,\gamma)$, $\alpha_2=\nobreak\alpha_2(\beta)$ and $\alpha_3=\alpha_3(\beta,\gamma)$ as follows:
\begin{equation} \label{AAA}
\alpha_1(\beta,\gamma):=-(1-\beta)\gamma, \quad\alpha_2(\beta):=-(1-\beta)N, \quad \alpha_3(\beta,\gamma):=-\dfrac{N+(1-2\beta)\gamma}{2}.
\end{equation}

\proof[Proof of Theorem \ref{teo4}]
For brevity, define
$$\Lambda_\infty:=\esssup_{r>R_2} \ \dfrac{K(r)}{r^{\alpha_\infty} V(r)^{\beta_\infty}} \quad \text{and} \quad \lambda_\infty:=\essinf_{r>R_2} \ r^{\gamma_\infty}V(r).$$
Take $u\in H^2_{V,r}$ and $h\in H^2_{V}$ such that $\left \| u \right \| = \left \| h \right \| =1$.
Let $R\geq R_2$. 
Hereafter $C$ denotes any positive constant independent from $u$, $h$ and $R$.
For all $\xi\geq0$ we have
\begin{equation} \label{leLam}
\begin{split}
\esssup_{r>R} \ \dfrac{K(r)}{r^{\alpha_\infty+\xi\gamma_\infty} V(r)^{\beta_\infty+\xi}} \leq \esssup_{r>R_2} \ \dfrac{K(r)}{r^{\alpha_\infty} V(r)^{\beta_\infty}(r^{\gamma_\infty}V(r))^\xi} 
 \leq \dfrac{\Lambda_\infty}{\lambda_\infty^\xi }\leq +\infty.
\end{split}
\end{equation}
We now distinguish several cases. In each of them we will choose a suitable  $\xi\geq0$ and we will apply Lemma \ref{lemma2} with $\Omega=B_R^c$, $\alpha=\alpha_\infty+\xi\gamma_\infty$, $\beta=\beta_\infty+\xi$, $m=c_\infty\lambda_\infty^{-1/4} \left \| u \right \|=c_\infty\lambda_\infty^{-1/4}$, $\nu=\frac{2(N-2)-\gamma_\infty}{4}$  and
$$\Lambda=\esssup_{r>R} \ \dfrac{K(r)}{r^{\alpha_\infty+\xi\gamma_\infty} V(r)^{\beta_\infty+\xi}}.$$
In each case we will get
$$\int_{B^c_{R}} K(|x|)|u|^{q_2-1}|h| \, dx \leq CR^\delta$$
for some $\delta<0$, independent from $R$, whence the thesis follows. 
Recalling the definitions \eqref{AAA}, we set $\alpha_1=\alpha_1(\beta_\infty,\gamma_\infty)$, $\alpha_2=\alpha_2(\beta_\infty)$ and $\alpha_3=\alpha_3(\beta_\infty,\gamma_\infty)$ for brevity.\smallskip

\noindent$\bullet$ Case $ \ \alpha_\infty\geq\alpha_1$. We take $\xi=1-\beta_\infty$ and apply Lemma \ref{lemma2} with $\beta=\beta_\infty+\xi=1$ and $\alpha=\alpha_\infty+\xi\gamma_\infty=\alpha_\infty+(1-\beta_\infty)\gamma_\infty$. We get
\begin{equation*}
\begin{split}
\int_{B^c_{R}} K(|x|)|u|^{q_2-1}|h| \, dx & \leq C\left(\int_{B^c_{R}} |x|^{2\alpha-2\nu(q_2-2)}V(|x|)|u|^{2} \, dx \right)^\frac{1}{2} \\
& \leq C\left(R^{2\alpha-2\nu(q_2-2)}\int_{B^c_{R}} V(|x|)|u|^{2} \, dx \right)^\frac{1}{2} 
\leq C R^{\alpha-\nu(q_2-2)},
\end{split}
\end{equation*}
because
\begin{equation*}
\begin{split}
\alpha-\nu(q_2-2)&=\alpha_\infty+(1-\beta_\infty)\gamma_\infty-\dfrac{2(N-2)-\gamma_\infty}{4}(q_2-2)\\
&=\dfrac{2\alpha_\infty+\gamma_\infty-2\beta_\infty\gamma_\infty+2N-4}{2}-\dfrac{2(N-2)-\gamma_\infty}{4}q_2 \\
&=\dfrac{2(N-2)-\gamma_\infty}{4}\left(2\dfrac{2\alpha_\infty+(1-2\beta_\infty)\gamma_\infty+2(N-2)}{2(N-2)-\gamma_\infty}-q_2 \right) \\
&=\dfrac{2(N-2)-\gamma_\infty}{4}\left(q_{**}-q_2 \right)<0.
\end{split}
\end{equation*}

\noindent$\bullet$ Case $ \ \max\{\alpha_2,\alpha_3\}<\alpha_\infty<\alpha_1$. We take $\xi=\frac{\alpha_\infty+(1-\beta_\infty)N}{N-\gamma_\infty}$ and apply Lemma \ref{lemma2} with $\beta=\beta_\infty+\xi$ and $\alpha=\alpha_\infty+\xi\gamma_\infty$. Notice that, if $\alpha_3<\alpha_\infty<\alpha_1$, then
$$\beta=\beta_\infty+\xi=\frac{\alpha_\infty-\gamma_\infty\beta_\infty+N}{N-\gamma_\infty}\in \left(\dfrac{1}{2},1\right).$$
On the other hand, if $\alpha_\infty=\alpha_2$ ($=\max\{\alpha_2,\alpha_3\}$ when $\frac{1}{2}<\beta_\infty<1$), then $\xi=0$ and
$$\beta=\beta_\infty\in \left(\dfrac{1}{2},1\right).$$
We obtain
\begin{equation*}
\begin{split}
\int_{B^c_{R}} K(|x|)|u|^{q_2-1}|h| \, dx  \leq C\left(\int_{B^c_{R}} |x|^\frac{\alpha-\nu(q_2-2\beta)}{1-\beta} \, dx \right)^{1-\beta} 
\leq C \left(R^{\frac{\alpha-\nu(q_2-2\beta)}{1-\beta}+N}\right)^{1-\beta},
\end{split}
\end{equation*}
because
\begin{equation*}
\begin{split}
\dfrac{\alpha-\nu(q_2-2\beta)}{1-\beta}+N =\dfrac{\nu}{1-\beta}\left(2\dfrac{\alpha_\infty-\beta_\infty\gamma_\infty+N}{N-\gamma_\infty}-q_2 \right) 
= \dfrac{\nu}{1-\beta}(q_*-q_2)<0.
\end{split}
\end{equation*}

\noindent$\bullet$ Case $\ \alpha_\infty\leq0=\alpha_2$ ($=\max\{\alpha_2,\alpha_3\}$) and $\beta_\infty=1$. We take $\xi=0$ and apply Lemma \ref{lemma2} with $\beta=\beta_\infty+\xi=1$ and $\alpha=\alpha_\infty+\xi\gamma_\infty=\alpha_\infty$. We get
\begin{equation*}
\begin{split}
\int_{B^c_{R}} K(|x|)|u|^{q_2-1}|h| \, dx  \leq C\left(\int_{B^c_{R}} |x|^{2\alpha_\infty-2\nu(q_2-2)}V(|x|)|u|^{2} \, dx \right)^\frac{1}{2} 
 \leq C R^{\alpha_\infty-\nu(q_2-2)},
\end{split}
\end{equation*}
because $\alpha_\infty-\nu(q_2-2)\leq-\nu(q_2-2)<0$. \smallskip

\noindent$\bullet$ Case $ \ \alpha_\infty\leq\alpha_2$ ($=\max\{\alpha_2,\alpha_3\}$) and $\frac{1}{2}<\beta_\infty<1$. We take $\xi=0$ and apply Lemma \ref{lemma2} with $\beta=\beta_\infty\in\left(\frac{1}{2},1\right)$ and $\alpha=\alpha_\infty+\xi\gamma_\infty=\alpha_\infty$. We get
\begin{equation*}
\begin{split}
\int_{B^c_{R}} K(|x|)|u|^{q_2-1}|h| \, dx  \leq C\left(\int_{B^c_{R}} |x|^\frac{\alpha_\infty-\nu(q_2-2\beta_\infty)}{1-\beta_\infty} \, dx \right)^{1-\beta_\infty} 
 \leq C \left(R^{\frac{\alpha_\infty-\nu(q_2-2\beta_\infty)}{1-\beta_\infty}+N}\right)^{1-\beta_\infty},
\end{split}
\end{equation*}
because
\begin{equation*}
\begin{split}
\dfrac{\alpha_\infty-\nu(q_2-2\beta_\infty)}{1-\beta_\infty}+N =\dfrac{\alpha_\infty+(1-\beta_\infty)N-\nu(q_2-2\beta_\infty)}{1-\beta_\infty} 
=\dfrac{\alpha_\infty-\alpha_2-\nu(q_2-2\beta_\infty)}{1-\beta_\infty}<0.
\end{split}
\end{equation*}

\noindent$\bullet$ Case $ \ \alpha_\infty\leq\alpha_3$ ($=\max\{\alpha_2,\alpha_3\}$) and $\beta_\infty\leq\frac{1}{2}$. We take $\xi=\frac{1-2\beta_\infty}{2}\geq0$ and apply Lemma \ref{lemma2} with $\beta=\beta_\infty+\xi=\frac{1}{2}$ and $\alpha=\alpha_\infty+\xi\gamma_\infty=\alpha_\infty$. We get
$$\int_{B^c_{R}} K(|x|)|u|^{q_2-1}|h| \, dx \leq C\left(\int_{B^c_{R}} |x|^{2\alpha-2\nu(q_2-1)} \, dx \right)^\frac{1}{2}\leq C R^{\alpha-\nu(q_2-1)+\frac{N}{2}}$$
because
$\alpha-\nu(q_2-1)+\frac{N}{2} =\alpha_\infty+\frac{1-2\beta_\infty}{2}\gamma_\infty+\frac{N}{2}-\nu(q_2-1) 
= \alpha_\infty-\alpha_3-\nu(q_2-1)<0.$
\endproof

\proof[Proof of Theorem \ref{teo5}]
Define
$$\Lambda_0:=\esssup_{r\in(0,R_1)} \ \dfrac{K(r)}{r^{\alpha_0} V(r)^{\beta_0}} \quad\text{and} \quad\lambda_0:=\essinf_{r\in(0,R_1)} \ r^{\gamma_0}V(r).$$
If $\gamma_0=4$ the thesis derives from Theorem \ref{teo2}, so we assume $\gamma>4$. To prove the result, we want to find a function $b(R)>0$ such that $b(R)\rightarrow0$ when $R\rightarrow0^+$ and
$$\int_{B_{R}} K(|x|)|u|^{q_1-1}|h| \, dx \leq b(R) \|u\|^{q_1-1}\|h\|, \ \ \forall u\in H^2_{V,r}, \ \forall h\in H^2_{V}.$$
So we fix $0<R\leq R_1$. Then
\begin{equation} \label{lambR}
\lambda(R):=\essinf_{r\in(0,R)} \ r^{\gamma_0}V(r)\geq \lambda_0>0
\end{equation}
and for all $\xi\geq0$ we have
\begin{equation} \label{Lambone}
\begin{split}
\Lambda_{\alpha_0+\xi\gamma_0,\beta_0+\xi}:=\esssup_{r\in(0,R)} \ \dfrac{K(r)}{r^{\alpha_0+\xi\gamma_0} V(r)^{\beta_0+\xi}} 
\leq \esssup_{r\in(0,R_1)} \ \dfrac{K(r)}{r^{\alpha_0} V(r)^{\beta_0}(r^{\gamma_0} V(r))^{\xi}}\leq \dfrac{\Lambda_0}{\lambda_0^\xi}<+\infty.
\end{split}
\end{equation}
We now consider several cases.\smallskip

\noindent$\bullet$ Case $ \ 4<\gamma_0<N$. In this case $(\alpha_0,q_1)\in \mathcal{A}_{\beta_0,\gamma_0}$ implies
$\alpha_0>\max\{\alpha_2,\alpha_3\}$ and
$$\max\{1,2\beta_0\}<q_1<\min\left\lbrace 2\dfrac{\alpha_0-\gamma_0\beta_0+N}{N-\gamma_0}, 2\dfrac{2\alpha_0+(1-2\beta_0)\gamma_0+2N-4}{2N-4-\gamma_0}\right\rbrace.$$
Hence, we can find $\xi\geq0$, independent from $R$, $u$ and $h$, such that $\alpha=\alpha_0+\xi\gamma_0$ and $\beta=\beta_0+\xi$ satisfy
\begin{equation} \label{condMagica}
\frac{1}{2}\leq\beta\leq1 \quad\text{and} \quad 2\beta<q_1<\dfrac{4\alpha+4N-2(\gamma_0+4)\beta}{2N-4-\gamma_0}.
\end{equation}
Recalling \eqref{lambR} and \eqref{Lambone}, we can apply Lemma \ref{lemma5} (with $q=q_1$), whence $\forall u\in H^2_{V,r}$ and $\forall h\in H^2_{V}$ we get
$$\int_{B_{R}} K(|x|)|u|^{q_1-1}|h| \, dx \leq c_0^{q_1-2\beta} a(R) R^\frac{4\alpha+4N-2(\gamma_0+4)\beta-(2N-4-\gamma_0)q_1}{4} \|u\|^{q_1-1}\|h\|.$$
This implies the thesis because $R^{4\alpha+4N-2(\gamma_0+4)\beta-(2N-4-\gamma_0)q_1}\rightarrow0$ as $R\rightarrow0^+$ and
\begin{equation*}
\begin{split}
a(R)=\Lambda_{\alpha_0+\xi\gamma_0,\beta_0+\xi}(R)\left( \dfrac{1}{\sqrt{\lambda(R)}}+\dfrac{R^\frac{\gamma_0-4}{2}}{\lambda(R)} \right)^\frac{q_1-2\beta}{2}
 \leq \dfrac{\Lambda_0}{\lambda_0^\xi}\left( \dfrac{1}{\sqrt{\lambda_0}}+\dfrac{R_1^\frac{\gamma_0-4}{2}}{\lambda_0} \right)^\frac{q_1-2\beta}{2}.
\end{split}
\end{equation*}

\noindent$\bullet$ Case $ \ N\leq\gamma_0<2N-4$. Again, $(\alpha_0,q_1)\in \mathcal{A}_{\beta_0,\gamma_0}$ implies that we can find $\xi\geq0$ such that $\alpha=\alpha_0+\xi\gamma_0$ and $\beta=\beta_0+\xi$ satisfy \eqref{condMagica}. We get the thesis applying Lemma \ref{lemma5}.\smallskip

\noindent$\bullet$ Case $ \ \gamma_0=2N-4$. In this case, from $(\alpha_0,q_1)\in \mathcal{A}_{\beta_0,\gamma_0}$ we infer that there exists $\xi\geq0$ such that $\alpha=\alpha_0+\xi\gamma_0$ and $\beta=\beta_0+\xi$ satisfy
$$\frac{1}{2}\leq\beta\leq1, \quad q_1>2\beta \quad\text{and} \quad 0<2\alpha+2N-(\gamma_0+4)\beta .$$
As in the previous cases, the thesis follows from Lemma \ref{lemma5}. \smallskip

\noindent$\bullet$ Case $ \ \gamma_0>2N-4$. In this case, the hypothesis $(\alpha_0,q_1)\in \mathcal{A}_{\beta_0,\gamma_0}$ implies that we can find $\xi\geq0$ such that $\alpha=\alpha_0+\xi\gamma_0$ and $\beta=\beta_0+\xi$ satisfy
$$\frac{1}{2}\leq\beta\leq1 \quad\text{and} \quad q_1>\max\left\lbrace 2\beta,2\dfrac{2\alpha+2N-(\gamma_0+4)\beta}{2N-4-\gamma_0}\right\rbrace.$$
Again, the thesis follows from Lemma \ref{lemma5}.
\endproof

\section{Application to the bilaplacian equation \label{SEC: ex}}

In this section we state our existence results for Eq. \eqref{prob1}, which are Theorems \ref{THM: super} and \ref{THM: sub} below (see also Remark \ref{RMK: extend}). 
We let $N\geq 5$ and assume that $V$, $K$
and $Q$ satisfy $\left( \mathbf{V}\right) $, $\left( \mathbf{K}\right) $
with $s>\frac{2N}{N+4}$ (cf. Lemma \ref{LEM: symm crit} below) and the
following hypothesis:

\begin{itemize}
\item[$\left( \mathbf{Q}\right) $]  $Q:\mathbb{R}_{+}\rightarrow \left[
0,+\infty \right) $ is a measurable function such that the linear funtional $%
h\mapsto \int_{\mathbb{R}^{N}}Q\left( \left| x\right| \right) h\,dx$ is
continuous on $H_{V}^{2}.$
\end{itemize}

\noindent We also assume that $f:\mathbb{R}\rightarrow \mathbb{R}$ is a
continuous function satisfying the following condition, where $q_{1},q_{2}$
will be specified later:

\begin{itemize}
\item[$\left( f_{q_{1},q_{2}}\right) $]  $\exists M>0$ such that $\left|
f\left( t\right) \right| \leq M\min \left\{ t^{q_{1}-1},t^{q_{2}-1}\right\} $
for all $t\geq 0.$
\end{itemize}

\begin{rem} 
\begin{enumerate}
\item  Assumption $\left( \mathbf{Q}\right) $ is quite abstract, but it is
easy to find explicit conditions on $Q$ ensuring it. For example, by Rellich inequality, 
if $Q\in L^{2}(\mathbb{R}_{+},r^{N+3}dr)$ then one has 
\[
\left| \int_{\mathbb{R}^{N}}Q\left( \left| x\right| \right) h\,dx\right| \leq
\left( \int_{\mathbb{R}^{N}}Q\left( \left| x\right| \right) ^{2}\left| x\right|
^{4}dx\right) ^{1/2}\left( \int_{\mathbb{R}^{N}}\frac{\left| h\right| ^{2}}{%
\left| x\right| ^{4}}dx\right) ^{1/2}\leq \left( \mathrm{const.}\right)
\left\| h\right\| ,\quad \forall h\in H_{V}^{2}.
\]
In a similar way, $\left( \mathbf{Q}\right) $ holds true if $Q\in L^{2N/(N+4)}(\mathbb{R%
}_{+},r^{N-1}dr)$ (by Sobolev inequality) or $V^{-1/2}Q\in L^{2}(\mathbb{R}%
_{+},r^{N-1}dr)$ (by definition of $H_{V}^{2}$). Other similar conditions
ensuring the same result can be obtained by the interpolation Hardy-Sobolev
inequalities of \cite{hardy.sob.1, hardy.sob.2} (see also \cite{Pass-Ruf}).

\item  Assumption $\left(f_{q_{1},q_{2}}\right) $ implies $\left|
f\left( t\right) \right| \leq M\,t^{q-1}$ for all $t\geq 0$ and $q\in \left[
q_{1},q_{2}\right] $, whence it is more stringent than a single-power growth
assumption if $q_{1}\neq q_{2}$. On the other hand we will never require $%
q_{1}\neq q_{2}$, so that our results will also concern single-power
nonlinearities as long as we can take $q_{1}=q_{2}$ in $\left( f_{q_{1},q_{2}}\right) $.
\end{enumerate}
\end{rem}

We are interested in finding \emph{radial weak solutions} of Eq. \eqref{prob1},
i.e., functions $u\in H_{V,r}^{2}$ such that 
\begin{equation}
\int_{\mathbb{R}^{N}}\triangle u\cdot \triangle h\,dx+\int_{\mathbb{R}%
^{N}}V\left( \left| x\right| \right) uh\,dx=\int_{\mathbb{R}^{N}}K\left(
\left| x\right| \right) f\left( u\right) h\,dx+\int_{\mathbb{R}^{N}}Q\left(
\left| x\right| \right) h\,dx\quad \textrm{for all }h\in H_{V}^{2}.
\label{weak sol}
\end{equation}

Our existence results are the following.

\begin{thm}
\label{THM: super}Assume $Q=0$ and assume that there exist $q_{1},q_{2}>2$
such that $\left( f_{q_{1},q_{2}}\right) $ and $\left( \mathcal{R}%
_{q_{1},q_{2}}^{\prime \prime }\right) $ hold. Assume furthermore that $f$
satisfies:

\begin{itemize}
\item[$\left( f_{1}\right) $]  $\exists \theta >2$ such that $0\leq \theta
F\left( t\right) \leq f\left( t\right) t$ for all $t\geq 0;$

\item[$\left( f_{2}\right) $]  $\exists t_{0}>0$ such that $F\left(
t_{0}\right) >0.$
\end{itemize}

\noindent If $K\left( \left| \cdot \right| \right) \in L^{1}(\mathbb{R}^{N})$,
we can replace assumptions $\left( f_{1}\right) $-$\left( f_{2}\right) $
with the following one:

\begin{itemize}
\item[$\left( f_{3}\right) $]  $\exists \theta >2$ and $\exists t_{0}>0$
such that $0<\theta F\left( t\right) \leq f\left( t\right) t$ for all $t\geq
t_{0}.$
\end{itemize}

\noindent Then Eq. \eqref{prob1} has a nonzero nonnegative radial weak solution.
\end{thm}

\begin{thm}
\label{THM: sub}Assume that there exist $q_{1},q_{2}\in \left( 1,2\right) $
such that $\left( f_{q_{1},q_{2}}\right) $ and $\left( \mathcal{R}%
_{q_{1},q_{2}}^{\prime \prime }\right) $ hold. Assume furthermore that
either $Q\neq 0$ (meaning that $Q$ does not vanish almost everywhere), or $%
Q=0$ and $f$ satisfies the following condition:

\begin{itemize}
\item[$\left( f_{4}\right) $]  $\exists \theta <2$ and $\exists t_{0},m>0$
such that $F\left( t\right) \geq mt^{\theta }$ for all $0\leq t\leq t_{0}.$
\end{itemize}

\noindent If $Q\neq 0$, we also allow the case $\max \left\{
q_{1},q_{2}\right\} =2>\min \left\{ q_{1},q_{2}\right\} >1$. Then Eq. \eqref{prob1}
has a nonzero nonnegative radial weak solution.
\end{thm}

The above existence results will be proved in Section \ref{SEC: app pf} and
can be generalized and complemented by other results in different and quite
standard ways (see Remark \ref{RMK: extend} below). They rely on assumption $%
\left( \mathcal{R}_{q_{1},q_{2}}^{\prime \prime }\right) $, which is rather
abstract but, as already discussed in Section \ref{SEC:MAIN}, it can be
granted in concrete cases through Theorems \ref{teo2}-\ref{teo5}, which ensure $\left( 
\mathcal{R}_{q_{1},q_{2}}^{\prime \prime }\right) $ for suitable ranges of
exponents $q_{1}$ and $q_{2}$ by explicit conditions on the potentials. Some
basic examples of nonlinearities satisfying $\left( f_{q_{1},q_{2}}\right) $
and the other assumptions of our results can be found in \cite[Example 4.11]
{bad3}.

\begin{rem} 
\begin{enumerate}
\item  In Theorem \ref{THM: super}, the information $K\left( \left| \cdot
\right| \right) \in L^{1}(\mathbb{R}^{N})$ actually allows weaker hypotheses on
the nonlinearity, as assumptions $\left( f_{1}\right) $ and $\left(
f_{2}\right) $ imply $\left( f_{3}\right) $.

\item  In Theorem \ref{THM: sub}, the case $\max \left\{ q_{1},q_{2}\right\}
=2>\min \left\{ q_{1},q_{2}\right\} >1$ cannot be considered if $\left(
f_{4}\right) $ holds, as $\left( f_{4}\right) $ and $\left(
f_{q_{1},q_{2}}\right) $ imply $\max \left\{ q_{1},q_{2}\right\} \leq \theta
<2$.
\end{enumerate}
\end{rem}

\begin{rem}
\label{RMK: extend} 
\begin{enumerate}
\item  Theorems \ref{THM: super} and \ref{THM: sub} can be easily adapted to
the case of Eq. \eqref{prob1} with a general right hand term $g\left( \left|
x\right| ,u\right) $ (see \cite[Section 3]{bad2} and \cite[Section 4]
{bad3}). Moreover, they can be complemented with multiplicity results by
standard variational techniques (see again \cite[Section 3]{bad2} and 
\cite[Section 4]{bad3}).

\item  Theorems \ref{THM: super} and \ref{THM: sub} can be used to derive
existence results for Eq. \eqref{prob1} with Dirichlet boundary conditions
in bounded balls or exterior radial domains, by suitably modifying the
potentials $V$ and $K$ in order to reduce the Dirichlet problem to the
problem in $\mathbb{R}^{N}$. In this cases, a single-power growth condition on
the nonlinearity is sufficient and, respectively, only assumptions on the
potentials near the origin or at infinity are needed. We leave the details
to the interested reader, which we refer to \cite[Section 5]{bad2} for
similar results and related arguments.

\item  Using some ideas of \cite{demarque-miya}, we think that our compactness and existence results can be easily extended to the case of inhomogeneous bilaplacian equations of the form 
\[
\Delta^2 u+V(|x|)u^{p-1}=K(|x|)f(u)+Q(|x|) \quad\text{in }\mathbb{R}^N
\]
with $1<p<N$, $p\neq 2$.
\end{enumerate}
\end{rem}

\section{Proof of Theorems \ref{THM: super} and \ref{THM: sub}}  \label{SEC: app pf}

In this section we apply the compactness results of Section \ref{SEC:MAIN} 
to prove the existence results of Section \ref{SEC: ex}.
Let $N\geq 5$ and assume that $V$, $K$ and $Q$ satisfy $\left( \mathbf{V}%
\right) $, $\left( \mathbf{K}\right) $ and $\left( \mathbf{Q}\right) $. Let $%
f:\mathbb{R}\rightarrow \mathbb{R}$ be a continuous function and set $%
F\left( t\right) :=\int_{0}^{t}f\left( s\right) ds$.

The weak solutions of Eq. \eqref{prob1} are (at least formally) the critical points
of the functional 
\begin{equation}
I\left( u\right) :=\frac{1}{2}\left\| u\right\| ^{2}-\int_{\mathbb{R}%
^{N}}K\left( \left| x\right| \right) F\left( u\right) dx-\int_{\mathbb{R}%
^{N}}Q\left( \left| x\right| \right) u\,dx.  \label{I:=}
\end{equation}
As a matter of fact, by the continuous embedding of Theorem \ref{teo1}
and the results of \cite{bad4} about Nemytski\u{\i} operators on the sum of
Lebesgue spaces, (\ref{I:=}) defines a $C^{1}$ functional on $H_{V,r}^{2}$
provided that $\left( f_{q_{1},q_{2}}\right) $ and $\left( \mathcal{S}%
_{q_{1},q_{2}}^{\prime }\right) $ hold for some $q_{1},q_{2}>1$. In this
case, the Fr\'{e}chet derivative of $I$ at any $u\in H_{V,r}^{2}$ is given
by 
\begin{equation}
I^{\prime }\left( u\right) h=\int_{\mathbb{R}^{N}}\left( \triangle u\cdot
\triangle h\,+V\left( \left| x\right| \right) uh\right) dx-\int_{\mathbb{R}%
^{N}}\left( K\left( \left| x\right| \right) f\left( u\right) +Q\left( \left|
x\right| \right) \right) h\,dx\,,\quad \forall h\in H_{V,r}^{2}\,,
\label{I'(u)h=}
\end{equation}
but $I$ does not need to be well defined on the whole space $H_{V}^{2}$, and
therefore the classical Palais' Principle of Symmetric Criticality \cite
{Palais} does not actually ensure that the critical points of $%
I:H_{V,r}^{2}\rightarrow \mathbb{R}$ are weak solutions of Eq. \eqref{prob1}. This is
the aim of our first lemma, which relies on the following stronger version
of condition $\left( \mathcal{S}_{q_{1},q_{2}}^{\prime }\right) $:

\begin{itemize}
\item[$\left( \mathcal{R}_{q_{1},q_{2}}^{\prime }\right) $]  $\exists
R_{1},R_{2}>0$ such that $\mathcal{R}_{0}\left( q_{1},R_{1}\right) <\infty $
and $\mathcal{R}_{\infty }\left( q_{2},R_{2}\right) <\infty .$
\end{itemize}

\begin{lem}
\label{LEM: symm crit}Assume $s>\frac{2N}{N+4}$ in condition $(\mathbf{K})$
and assume that there exist $q_{1},q_{2}>1$ such that $\left(
f_{q_{1},q_{2}}\right) $ and $\left( \mathcal{R}_{q_{1},q_{2}}^{\prime
}\right) $ hold. Then every critical point of $I:H_{V,r}^{2}\rightarrow \mathbb{%
R}$ is a weak solution to Eq. \eqref{prob1}.
\end{lem}

\proof
Let $u\in H_{V,r}^{2}$ and assume $R_{1}<R_{2}$ in $\left( \mathcal{R}%
_{q_{1},q_{2}}^{\prime }\right) $, which is not restrictive by the
monotonicity of $\mathcal{R}_{0}$ and $\mathcal{R}_{\infty }$. By Lemma \ref{lemma1}, 
there exists a constant $C>0$ (also dependent on $u$) such
that for all $h\in H_{V}^{2}$ we have 
\[
\int_{B_{R_{2}}\setminus B_{R_{1}}}K\left( \left| x\right| \right) \left|
u\right| ^{q_{1}-1}\left| h\right| dx\leq C\left\| h\right\| 
\]
and therefore, by $\left( f_{q_{1},q_{2}}\right) $, we get 
\begin{eqnarray*}
\int_{\mathbb{R}^{N}}K\left( \left| x\right| \right) \left| f\left( u\right)
\right| \left| h\right| dx &\leq &M\int_{\mathbb{R}^{N}}K\left( \left| x\right|
\right) \min \{\left| u\right| ^{q_{1}-1},\left| u\right| ^{q_{2}-1}\}\left|
h\right| dx \\
&\leq &M\left( \int_{B_{R_{1}}}K\left( \left| x\right| \right) \left|
u\right| ^{q_{1}-1}\left| h\right| dx+\int_{B_{R_{2}}^{c}}K\left( \left|
x\right| \right) \left| u\right| ^{q_{2}-1}\left| h\right| dx+C\left\|
h\right\| \right) \\
&\leq &M\left( \left\| u\right\| ^{q_{1}-1}\mathcal{R}_{0}\left(
q_{1},R_{1}\right) +\left\| u\right\| ^{q_{2}-1}\mathcal{R}_{\infty }\left(
q_{2},R_{2}\right) +C\right) \left\| h\right\| .
\end{eqnarray*}
Together with assumption $\left( \mathbf{Q}\right) $, this gives that the
linear operator 
\[
T\left( u\right) h:=\int_{\mathbb{R}^{N}}\left( \triangle u\cdot \triangle
h\,+V\left( \left| x\right| \right) uh\right) dx-\int_{\mathbb{R}^{N}}\left(
K\left( \left| x\right| \right) f\left( u\right) +Q\left( \left| x\right|
\right) \right) h\,dx 
\]
is well defined and continuous on $H_{V}^{2}$. Hence, by Riesz
representation theorem, there exists a unique $\tilde{u}\in H_{V}^{2}$ such
that $T\left( u\right) h=\left( \tilde{u},h\right) $ for all $h\in H_{V}^{2}$%
, where $\left( \cdot ,\cdot \right) $ is the scalar product defined in
(\ref{inn.prod}). By means of obvious changes of variables it is easy to infer that $%
\tilde{u}\in H_{V,r}^{2}$, so that $T\left( u\right) =0$ on $H_{V,r}^{2}$
implies $\tilde{u}=0$ and hence (\ref{weak sol}).
\endproof
\bigskip

Hereafter, we will assume that the hypotheses of Theorems \ref{THM: super}
and \ref{THM: sub} also include the following assumptions respectively: $%
f\left( t\right) =0$ for $t<0$ in Theorem \ref{THM: super}, and $f$ is odd
in Theorem \ref{THM: sub}. This can be done without restriction, since
Theorems \ref{THM: super} and \ref{THM: sub} concern nonnegative solutions
and all their assumptions still hold true if we replace $f\left( t\right) $
respectively with $f\left( t\right) \chi _{\mathbb{R}_{+}}\left( t\right) $ and 
$f\left( t\right) \chi _{\mathbb{R}_{+}}\left( t\right) -f\left( \left|
t\right| \right) \chi _{\mathbb{R}_{-}}\left( t\right) $ ($\chi _{\mathbb{R}_{\pm
}}$ denotes the characteristic function of $\mathbb{R}_{\pm }$).

With such additional assumptions, $\left( f_{q_{1},q_{2}}\right) $ implies
that there exists $\tilde{M}>0$ such that 
\begin{equation}
\left| F\left( t\right) \right| \leq \tilde{M}\min \left\{ \left| t\right|
^{q_{1}},\left| t\right| ^{q_{2}}\right\} \quad \text{for all }t\in \mathbb{R}.
\label{G_pq}
\end{equation}

\begin{lem}
\label{LEM:MPgeom}Let $L_{0}$ be the norm of the linear continuous
functional $h\in H_{V}^{2}\mapsto \int_{\mathbb{R}^{N}}Q\left( \left| x\right|
\right) h\,dx$. If $\left( f_{q_{1},q_{2}}\right) $ and $\left( \mathcal{S}%
_{q_{1},q_{2}}^{\prime }\right) $ hold for some $q_{1},q_{2}>1$, then there
exist two constants $c_{1},c_{2}>0$ such that 
\begin{equation}
I\left( u\right) \geq \frac{1}{2}\left\| u\right\| ^{2}-c_{1}\left\|
u\right\| ^{q_{1}}-c_{2}\left\| u\right\| ^{q_{2}}-L_{0}\left\| u\right\|
\qquad \text{for all }u\in H_{V,r}^{2}.  \label{LEM:MPgeom: th}
\end{equation}
If $\left( \mathcal{S}_{q_{1},q_{2}}^{\prime \prime }\right) $ also holds,
then $\forall \varepsilon >0$ there exist two constants $c_{1}\left(
\varepsilon \right) ,c_{2}\left( \varepsilon \right) >0$ such that (\ref
{LEM:MPgeom: th}) holds both with $c_{1}=\varepsilon $, $c_{2}=c_{2}\left(
\varepsilon \right) $ and with $c_{1}=c_{1}\left( \varepsilon \right) $, $%
c_{2}=\varepsilon $.
\end{lem}

\proof
Let $i\in \left\{ 1,2\right\} $ and assume $R_{1}<R_{2}$ in $\left( \mathcal{%
S}_{q_{1},q_{2}}^{\prime }\right) $, which is not restrictive by the
monotonicity of $\mathcal{S}_{0}$ and $\mathcal{S}_{\infty }$. By Lemma \ref{lemma1} 
and the continuous embedding $H_{V,r}^{2}\hookrightarrow L_{%
\mathrm{loc}}^{q_{i}}(\mathbb{R}^{N})$, there exists a constant $%
c_{R_{1},R_{2}}^{\left( i\right) }>0$ such that for all $u\in H_{V,r}^{2}$
we have 
\[
\int_{B_{R_{2}}\setminus B_{R_{1}}}K\left( \left| x\right| \right) \left|
u\right| ^{q_{i}}dx\leq c_{R_{1},R_{2}}^{\left( i\right) }\left\| u\right\|
^{q_{i}} 
\]
and therefore, by (\ref{G_pq}), 
\begin{eqnarray}
&&\left| \int_{\mathbb{R}^{N}}\left( K\left( \left| x\right| \right) F\left(
u\right) +Q\left( \left| x\right| \right) u\right) dx\right|  \nonumber \\
&\leq &\int_{\mathbb{R}^{N}}K\left( \left| x\right| \right) \left| F\left(
u\right) \right| dx+\left| \int_{\mathbb{R}^{N}}Q\left( \left| x\right| \right)
u\,dx\right| \leq \tilde{M}\int_{\mathbb{R}^{N}}K\left( \left| x\right| \right)
\min \left\{ \left| u\right| ^{q_{1}},\left| u\right| ^{q_{2}}\right\}
dx+L_{0}\left\| u\right\|  \nonumber \\
&\leq &\tilde{M}\left( \int_{B_{R_{1}}}K\left( \left| x\right| \right)
\left| u\right| ^{q_{1}}dx+\int_{B_{R_{2}}^{c}}K\left( \left| x\right|
\right) \left| u\right| ^{q_{2}}dx+\int_{B_{R_{2}}\setminus
B_{R_{1}}}K\left( \left| x\right| \right) \left| u\right| ^{q_{i}}dx\right)
+L_{0}\left\| u\right\|  \nonumber \\
&\leq &\tilde{M}\left( \left\| u\right\| ^{q_{1}}\mathcal{S}_{0}\left(
q_{1},R_{1}\right) +\left\| u\right\| ^{q_{2}}\mathcal{S}_{\infty }\left(
q_{2},R_{2}\right) +c_{R_{1},R_{2}}^{\left( i\right) }\left\| u\right\|
^{q_{i}}\right) +L_{0}\left\| u\right\|  \label{LEM:MPgeom: pf} \\
&=&c_{1}\left\| u\right\| ^{q_{1}}+c_{2}\left\| u\right\|
^{q_{2}}+L_{0}\left\| u\right\| ,  \nonumber
\end{eqnarray}
with obvious definition of the constants $c_{1}$ and $c_{2}$, independent of 
$u$. This proves (\ref{LEM:MPgeom: th}). If $\left( \mathcal{S}%
_{q_{1},q_{2}}^{\prime \prime }\right) $ also holds, then $\forall
\varepsilon >0$ we can fix $R_{1,\varepsilon }<R_{2,\varepsilon }$ such that 
$\tilde{M}\mathcal{S}_{0}\left( q_{1},R_{1,\varepsilon }\right) <\varepsilon 
$ and $\tilde{M}\mathcal{S}_{\infty }\left( q_{2},R_{2,\varepsilon }\right)
<\varepsilon $, so that inequality (\ref{LEM:MPgeom: pf}) becomes 
\[
\left| \int_{\mathbb{R}^{N}}\left( K\left( \left| x\right| \right) F\left(
u\right) +Q\left( \left| x\right| \right) u\right) dx\right| \leq
\varepsilon \left\| u\right\| ^{q_{1}}+\varepsilon \left\| u\right\|
^{q_{2}}+c_{R_{1,\varepsilon },R_{2,\varepsilon }}^{\left( i\right) }\left\|
u\right\| ^{q_{i}}+L_{0}\left\| u\right\| . 
\]
The conclusion thus ensues by taking $i=1$ and $c_{1}\left( \varepsilon
\right) =\varepsilon +c_{R_{1,\varepsilon },R_{2,\varepsilon }}^{\left(
1\right) }$, or $i=2$ and $c_{2}\left( \varepsilon \right) =\varepsilon
+c_{R_{1,\varepsilon },R_{2,\varepsilon }}^{\left( 2\right) }$.%
\endproof

\begin{lem}
\label{LEM:PS}Under the assumptions of Theorem \ref{THM: super}, the
functional $I:H_{V,r}^{2}\rightarrow \mathbb{R}$ satisfies the Palais-Smale
condition.
\end{lem}

\proof
By $\left( f_{1}\right) $ or $\left( f_{3}\right) $ together with the
additional assumption $f\left( t\right) =0$ for $t<0$, we have that either $%
\left( f_{1}\right) $ holds for all $t\in \mathbb{R}$, or $K\left( \left| \cdot
\right| \right) \in L^{1}(\mathbb{R}^{N})$ and $f$ satisfies 
\begin{equation}
\theta F\left( t\right) \leq f\left( t\right) t\quad \text{for all~}\left|
t\right| \geq t_{0}.  \label{LEM:PS: AR}
\end{equation}
Let $\left\{ u_{n}\right\} $ be a sequence in $H_{V,r}^{2}$ such that $%
\left\{ I\left( u_{n}\right) \right\} $ is bounded and $I^{\prime }\left(
u_{n}\right) \rightarrow 0$ in the dual space of $H_{V,r}^{2}$. Then 
\[
\frac{1}{2}\left\| u_{n}\right\| ^{2}-\int_{\mathbb{R}^{N}}K\left( \left|
x\right| \right) F\left( u_{n}\right) dx=I\left( u_{n}\right) +\int_{\mathbb{R}%
^{N}}Q\left( \left| x\right| \right) u_{n}dx=O\left( 1\right) +O\left(
1\right) \left\| u_{n}\right\| 
\]
and 
\[
\left\| u_{n}\right\| ^{2}-\int_{\mathbb{R}^{N}}K\left( \left| x\right| \right)
f\left( u_{n}\right) u_{n}dx=I^{\prime }\left( u_{n}\right) u_{n}+\int_{\mathbb{%
R}^{N}}Q\left( \left| x\right| \right) u_{n}dx=o\left( 1\right) \left\|
u_{n}\right\| +O\left( 1\right) \left\| u_{n}\right\| . 
\]
If $f$ satisfies $\left( f_{1}\right) $, we get 
\begin{eqnarray*}
\frac{1}{2}\left\| u_{n}\right\| ^{2}+O\left( 1\right) +O\left( 1\right)
\left\| u_{n}\right\| &=&\int_{\mathbb{R}^{N}}K\left( \left| x\right| \right)
F\left( u_{n}\right) dx \\
&\leq &\frac{1}{\theta }\int_{\mathbb{R}^{N}}K\left( \left| x\right| \right)
f\left( u_{n}\right) u_{n}dx=\frac{1}{\theta }\left\| u_{n}\right\|
^{2}+o\left( 1\right) \left\| u_{n}\right\| +O\left( 1\right) \left\|
u_{n}\right\| ,
\end{eqnarray*}
which implies that $\left\{ \left\| u_{n}\right\| \right\} $ is bounded
since $\theta >2$. If $K\left( \left| \cdot \right| \right) \in L^{1}(\mathbb{R}%
^{N})$ and $f$ satisfies (\ref{LEM:PS: AR}), we have 
\begin{eqnarray*}
\int_{\left\{ \left| u_{n}\right| \geq t_{0}\right\} }K\left( \left|
x\right| \right) f\left( u_{n}\right) u_{n}dx &\leq &\int_{\mathbb{R}%
^{N}}K\left( \left| x\right| \right) f\left( u_{n}\right)
u_{n}dx+\int_{\left\{ \left| u_{n}\right| <t_{0}\right\} }K\left( \left|
x\right| \right) \left| f\left( u_{n}\right) u_{n}\right| dx \\
&\leq &\int_{\mathbb{R}^{N}}K\left( \left| x\right| \right) f\left(
u_{n}\right) u_{n}dx+M\int_{\left\{ \left| u_{n}\right| <t_{0}\right\}
}K\left( \left| x\right| \right) \min \left\{ \left| u_{n}\right|
^{q_{1}},\left| u_{n}\right| ^{q_{2}}\right\} dx \\
&\leq &\int_{\mathbb{R}^{N}}K\left( \left| x\right| \right) f\left(
u_{n}\right) u_{n}dx+M\min \left\{ t_{0}^{q_{1}},t_{0}^{q_{2}}\right\}
\left\| K\right\| _{L^{1}(\mathbb{R}^{N})},
\end{eqnarray*}
and then, using (\ref{G_pq}), we get 
\begin{eqnarray*}
&&\frac{1}{2}\left\| u_{n}\right\| ^{2}+O\left( 1\right) +O\left( 1\right)
\left\| u_{n}\right\| \\
&=&\int_{\left\{ \left| u_{n}\right| <t_{0}\right\} }K\left( \left| x\right|
\right) F\left( u_{n}\right) dx+\int_{\left\{ \left| u_{n}\right| \geq
t_{0}\right\} }K\left( \left| x\right| \right) F\left( u_{n}\right) dx \\
&\leq &\tilde{M}\int_{\left\{ \left| u_{n}\right| <t_{0}\right\} }K\left(
\left| x\right| \right) \min \left\{ \left| u_{n}\right| ^{q_{1}},\left|
u_{n}\right| ^{q_{2}}\right\} dx+\frac{1}{\theta }\int_{\left\{ \left|
u_{n}\right| \geq t_{0}\right\} }K\left( \left| x\right| \right) f\left(
u_{n}\right) u_{n}dx \\
&\leq &\tilde{M}\min \left\{ t_{0}^{q_{1}},t_{0}^{q_{2}}\right\} \left\|
K\right\| _{L^{1}(\mathbb{R}^{N})}+\frac{1}{\theta }\int_{\mathbb{R}^{N}}K\left(
\left| x\right| \right) f\left( u_{n}\right) u_{n}dx+\frac{M}{\theta }\min
\left\{ t_{0}^{q_{1}},t_{0}^{q_{2}}\right\} \left\| K\right\| _{L^{1}(\mathbb{R}%
^{N})} \\
&=&\left( \tilde{M}+\frac{M}{\theta }\right) \min \left\{
t_{0}^{q_{1}},t_{0}^{q_{2}}\right\} \left\| K\right\| _{L^{1}(\mathbb{R}^{N})}+%
\frac{1}{\theta }\left\| u_{n}\right\| ^{2}+o\left( 1\right) \left\|
u_{n}\right\| +O\left( 1\right) \left\| u_{n}\right\| .
\end{eqnarray*}
This yields again that $\left\{ \left\| u_{n}\right\| \right\} $ is bounded.
Now, thanks to assumption $\left( \mathbf{Q}\right) $ and since the
embedding $H_{V,r}^{2}\hookrightarrow L_{K}^{q_{1}}+L_{K}^{q_{2}}$ is
compact by Theorem \ref{teo1} and the functional $u\mapsto \int_{\mathbb{R}%
^{N}}K\left( \left| x\right| \right) F\left( u\right) dx$ is $C^{1}$ on $%
L_{K}^{q_{1}}+L_{K}^{q_{2}}$ by \cite[Proposition 3.8]{bad4}, it is a
standard exercise to conclude that $\left\{ u_{n}\right\} $ has a strongly
convergent subsequence in $H_{V,r}^{2}$.%
\endproof
\bigskip

\noindent \textbf{Proof of Theorem \ref{THM: super}.}\quad We mean to apply
the Mountain-Pass Theorem. To this end, from (\ref{LEM:MPgeom: th}) of Lemma 
\ref{LEM:MPgeom}, where $L_{0}=0$ and $q_{1},q_{2}>2$, we readily infer that 
$\exists \rho >0$ such that 
\begin{equation}
\inf_{u\in H_{V,r}^{2},\,\left\| u\right\| =\rho }I\left( u\right)
>0=I\left( 0\right) .  \label{mp-geom}
\end{equation}
Now we check that $\exists \bar{u}\in W_{r}$ such that $\left\| \bar{u}%
\right\| >\rho $ and $I\left( \bar{u}\right) <0$. To this end, from $\left(
f_{3}\right) $ (which holds in any case, as $\left( f_{1}\right) $ and $%
\left( f_{2}\right) $ imply $\left( f_{3}\right) $), we deduce that $F\left(
t\right) \geq t_{0}^{-\theta }F\left( t_{0}\right) t^{\theta }$ for all $%
t\geq t_{0}$. Then, taking into account assumption $\left( \mathbf{V}\right) 
$, we fix a nonnegative function $u_{0}\in C_{c}^{\infty }(\mathbb{R}^{N})\cap
H_{V,r}^{2}$ such that the set $\{x\in \mathbb{R}^{N}:u_{0}\left( x\right) \geq
t_{0}\}$ has positive Lebesgue measure. We now distinguish the case of
assumptions $\left( f_{1}\right) $ and $\left( f_{2}\right) $ from the case
with $K\left( \left| \cdot \right| \right) \in L^{1}(\mathbb{R}^{N})$. In the
first one, for every $\lambda >1$ we have 
\begin{eqnarray*}
\int_{\mathbb{R}^{N}}K\left( \left| x\right| \right) F\left( \lambda
u_{0}\right) dx &\geq &\int_{\left\{ \lambda u_{0}\geq t_{0}\right\}
}K\left( \left| x\right| \right) F\left( \lambda u_{0}\right) dx\geq \lambda
^{\theta }F\left( t_{0}\right) \int_{\left\{ \lambda u_{0}\geq t_{0}\right\}
}t_{0}^{-\theta }u_{0}^{\theta }dx \\
&\geq &\lambda ^{\theta }F\left( t_{0}\right) \int_{\left\{ u_{0}\geq
t_{0}\right\} }t_{0}^{-\theta }u_{0}^{\theta }dx\geq \lambda ^{\theta
}F\left( t_{0}\right) \int_{\left\{ u_{0}\geq t_{0}\right\} }dx>0
\end{eqnarray*}
and therefore, since $\theta >2$, we get 
\[
\lim_{\lambda \rightarrow +\infty }I\left( \lambda u_{0}\right) \leq
\lim_{\lambda \rightarrow +\infty }\left( \frac{\lambda ^{2}}{2}\left\|
u_{0}\right\| ^{2}-\lambda ^{\theta }F\left( t_{0}\right) \int_{\left\{
u_{0}\geq t_{0}\right\} }dx-\lambda \int_{\mathbb{R}^{N}}Q\left( \left|
x\right| \right) u_{0}dx\right) =-\infty . 
\]
If $K\left( \left| \cdot \right| \right) \in L^{1}(\mathbb{R}^{N})$, we observe
that (\ref{G_pq}) implies $F\left( t\right) \geq -\tilde{M}\min \left\{
t_{0}^{q_{1}},t_{0}^{q_{2}}\right\} $ for all $0\leq t\leq t_{0}$, so that,
arguing as above about the integral over $\left\{ \lambda u_{0}\geq
t_{0}\right\} $, for every $\lambda >1$ we obtain 
\begin{eqnarray*}
\int_{\mathbb{R}^{N}}K\left( \left| x\right| \right) F\left( \lambda
u_{0}\right) dx &=&\int_{\left\{ \lambda u_{0}<t_{0}\right\} }K\left( \left|
x\right| \right) F\left( \lambda u_{0}\right) dx+\int_{\left\{ \lambda
u_{0}\geq t_{0}\right\} }K\left( \left| x\right| \right) F\left( \lambda
u_{0}\right) dx \\
&\geq &-\tilde{M}\min \left\{ t_{0}^{q_{1}},t_{0}^{q_{2}}\right\}
\int_{\left\{ \lambda u_{0}<t_{0}\right\} }K\left( \left| x\right| \right)
dx+\lambda ^{\theta }F\left( t_{0}\right) \int_{\left\{ u_{0}\geq
t_{0}\right\} }dx,
\end{eqnarray*}
which implies 
\[
I\left( \lambda u_{0}\right) \leq \frac{\lambda ^{2}}{2}\left\|
u_{0}\right\| ^{2}+\tilde{M}\min \left\{ t_{0}^{q_{1}},t_{0}^{q_{2}}\right\}
\left\| K\right\| _{L^{1}(\mathbb{R}^{N})}-\lambda ^{\theta }F\left(
t_{0}\right) \int_{\left\{ u_{0}\geq t_{0}\right\} }dx-\lambda \int_{\mathbb{R}%
^{N}}Q\left( \left| x\right| \right) u_{0}dx\rightarrow -\infty 
\]
as $\lambda \rightarrow +\infty $. So, in any case, we can take $\bar{u}%
=\lambda u_{0}$ with $\lambda $ sufficiently large. As a conclusion, taking
into account Lemma \ref{LEM:PS}, the Mountain-Pass Theorem provides the
existence of a nonzero critical point $u\in H_{V,r}^{2}$ for $I$, which is a
weak solutions to Eq. \eqref{prob1} by Lemma \ref{LEM: symm crit}. Since the
additional assumption $f\left( t\right) =0$ for $t<0$ implies $I^{\prime
}\left( u\right) u_{-}=-\left\| u_{-}\right\| ^{2}$ (where $u_{-}\in
H_{V,r}^{2}$ is the negative part of $u$), one has $u_{-}=0$ and thus $u$ is
nonnegative.%
\endproof

\begin{lem}
\label{LEM:bdd}Under the assumptions of Theorem \ref{THM: sub}, the
functional $I:H_{V,r}^{2}\rightarrow \mathbb{R}$ is bounded from below and
coercive. In particular, if $Q=0$ and $f$ satisfies $\left( f_{4}\right) $,
then 
\begin{equation}
\inf_{v\in H_{V,r}^{2}}I\left( v\right) <0.  \label{LEM:bdd: inf<0}
\end{equation}
\end{lem}

\proof
$I$ is bounded below and coercive on $H_{V,r}^{2}$ thanks to Lemma \ref
{LEM:MPgeom}. Indeed, the result readily follows from (\ref{LEM:MPgeom: th})
if $q_{1},q_{2}\in \left( 1,2\right) $, while, if $\max \left\{
q_{1},q_{2}\right\} =2>\min \left\{ q_{1},q_{2}\right\} >1$, we fix $%
\varepsilon <1/2$ and use the second part of the lemma in order to get 
\[
I\left( u\right) \geq \left( \frac{1}{2}-\varepsilon \right) \left\|
u\right\| ^{2}-c\left( \varepsilon \right) \left\| u\right\| ^{\min \left\{
q_{1},q_{2}\right\} }-L_{0}\left\| u\right\| \qquad \text{for all }u\in
H_{V,r}^{2}, 
\]
which yields again the conclusion. In order to prove (\ref{LEM:bdd: inf<0})
if $Q=0$ and $\left( f_{4}\right) $ holds, we use assumption $\left( \mathbf{%
V}\right) $ to fix a nonzero function $u_{0}\in C_{c}^{\infty }(\mathbb{R}%
^{N})\cap H_{V,r}^{2}$ such that $0\leq u_{0}\leq t_{0}$. Then, by $\left(
f_{4}\right) $, for every $0<\lambda <1$ we get that $\lambda u_{0}\in
H_{V,r}^{2}$ satisfies 
\[
I\left( \lambda u_{0}\right) =\frac{1}{2}\left\| \lambda u_{0}\right\|
^{2}-\int_{\mathbb{R}^{N}}K\left( \left| x\right| \right) F\left( \lambda
u_{0}\right) dx\leq \frac{\lambda ^{2}}{2}\left\| u_{0}\right\| ^{2}-\lambda
^{\theta }m\int_{\mathbb{R}^{N}}K\left( \left| x\right| \right) u_{0}^{\theta
}dx. 
\]
Since $\theta <2$, this implies $I\left( \lambda u_{0}\right) <0$ for $%
\lambda $ sufficiently small and therefore (\ref{LEM:bdd: inf<0}) ensues.%
\endproof
\bigskip

\noindent \textbf{Proof of Theorem \ref{THM: sub}.}\quad 
Recall Lemma \ref{LEM:bdd} and let $\left\{ v_{n}\right\} $ be any minimizing
sequence for $\mu :=\inf_{v\in H_{V,r}^{2}}I\left( v\right) \in \mathbb{R}$. As 
$F$ is even and $Q$ is nonnegative, we have 
\begin{eqnarray}
I\left( \left| v_{n}\right| \right) &=&\frac{1}{2}\left\| v_{n}\right\|
^{2}-\int_{\mathbb{R}^{N}}K\left( \left| x\right| \right) F\left( \left|
v_{n}\right| \right) dx-\int_{\mathbb{R}^{N}}Q\left( \left| x\right| \right)
\left| v_{n}\right| dx \\
&=&\frac{1}{2}\left\| v_{n}\right\| ^{2}-\int_{\mathbb{R}^{N}}K\left( \left|
x\right| \right) F\left( v_{n}\right) dx-\int_{\left\{ v_{n}\geq 0\right\}
}Q\left( \left| x\right| \right) v_{n}dx+\int_{\left\{ v_{n}<0\right\}
}Q\left( \left| x\right| \right) v_{n}dx  \nonumber \\
&=&I\left( v_{n}\right) +2\int_{\left\{ v_{n}<0\right\} }Q\left( \left|
x\right| \right) v_{n}dx\leq I\left( v_{n}\right) ,  \nonumber
\end{eqnarray}
so that $\left| v_{n}\right| \in H_{V,r}^{2}$ is still a minimizing
sequence. Hence we can assume $v_{n}\geq 0$. Since $\left\{ v_{n}\right\} $
is bounded in $H_{V,r}^{2}$ by Lemma \ref{LEM:bdd} and the embedding $%
H_{V,r}^{2}\hookrightarrow L_{K}^{q_{1}}+L_{K}^{q_{2}}$ is compact by
assumption $\left( \mathcal{S}_{q_{1},q_{2}}^{\prime \prime }\right) $ and
Theorem \ref{teo1}, we can assume that there exists $u\in H_{V,r}^{2}$
such that (up to a subsequence) $v_{n}\rightharpoonup u$ in $H_{V,r}^{2}$
and $v_{n}\rightarrow u$ in $L_{K}^{q_{1}}+L_{K}^{q_{2}}$ and almost
everywhere. Hence $u$ is nonnegative and, thanks to $\left( \mathbf{Q}%
\right) $ and the continuity of the functional $v\mapsto \int_{\mathbb{R}%
^{N}}K\left( \left| x\right| \right) F\left( v\right) dx$ on $%
L_{K}^{q_{1}}+L_{K}^{q_{2}}$ (which follows from $\left(
f_{q_{1},q_{2}}\right) $ and \cite[Proposition 3.8]{bad4}), we have 
\[
\int_{\mathbb{R}^{N}}K\left( \left| x\right| \right) F\left( v_{n}\right)
dx+\int_{\mathbb{R}^{N}}Q\left( \left| x\right| \right) v_{n}dx\rightarrow
\int_{\mathbb{R}^{N}}K\left( \left| x\right| \right) F\left( u\right) dx+\int_{%
\mathbb{R}^{N}}Q\left( \left| x\right| \right) u\,dx. 
\]
By the weak lower semi-continuity of the norm, this implies 
\[
I\left( u\right) \leq \lim_{n\rightarrow \infty }\left( \frac{1}{2}\left\|
v_{n}\right\| ^{2}-\int_{\mathbb{R}^{N}}K\left( \left| x\right| \right) F\left(
v_{n}\right) dx+\int_{\mathbb{R}^{N}}Q\left( \left| x\right| \right)
v_{n}dx\right) =\mu 
\]
and thus $I\left( u\right) =\mu $. So $u$ is a critical point for $I$ and
thus a weak solutions to Eq. \eqref{prob1} by Lemma \ref{LEM: symm crit}. It remains
to show that $u\neq 0$. This is obvious if $Q=0$ and $f$ satisfies $\left(
f_{4}\right) $, since $\mu <0$ by Lemma \ref{LEM:bdd}. If $Q\neq 0$, assume
by contradiction that $u=0$. From (\ref{I'(u)h=}) we get $\int_{\mathbb{R}%
^{N}}Q\left( \left| x\right| \right) h\,dx=0$ for all $h\in H_{V,r}^{2}$ and
therefore $Q=0$, which is a contradiction.%
\endproof

\end{document}